\documentclass[12pt,reqno]{article}

\usepackage[usenames]{color}
\usepackage{graphicx}
\usepackage{xcolor}
\usepackage{fullpage}
\usepackage{amsmath,amssymb}
\usepackage{amsthm}
\usepackage{amsfonts}
\usepackage{latexsym}
\usepackage{epigraph}
\usepackage[colorlinks=true,
linkcolor=webgreen,
filecolor=webbrown,
citecolor=webgreen]{hyperref}

\definecolor{webgreen}{rgb}{0,.5,0}
\definecolor{webbrown}{rgb}{.6,0,0}
\usepackage{tikz}
\usetikzlibrary{matrix}
\usepackage{amsmath}
\usepackage{caption}

\usepackage{amscd}
\usepackage{float}
\usepackage{graphics}
\usepackage{epsf}

\newcommand{\seqnum}[1]{\href{https://oeis.org/#1}{\underline{#1}}}

\newcommand{\Ei}{{E_{\infty}}} 
\newcommand{\Ri}{{R_{\infty}}} 
\newcommand{\sB}{{\mathcal B}} 
\newcommand{\sC}{{\mathcal C}} 
\newcommand{\sG}{{\mathcal G}} 
\newcommand{\beql}[1]{\begin{equation}\label{#1}}
\newcommand{\eeq}{\end{equation}}

\newcommand{\cpa}{{\mathnormal cpa}}

\newcommand{\mabs}{{family}}
\newcommand{\mabp}{{families}}
\newcommand{\affineGroup}{\sG_{\mathrm{aff}}}
\newcommand{\simGroup}{\sG_{\mathrm{sim}}}
\newcommand{\isoGroup}{\sG_{\mathrm{iso}}}
\newcommand{\wu}{\includegraphics[height=1.6ex]{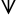}}
\newcommand{\qoppaletter}{
	\raisebox{-0.6ex}{\includegraphics[height=1.8ex]{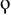}}
}
\newcommand{\runewynn}{
	\raisebox{-0.2ex}{\includegraphics[height=1.7ex]{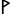}}
}

\usetikzlibrary{calc}
\usetikzlibrary{arrows.meta}
\usetikzlibrary{shapes.geometric}
\tikzset{axis line style/.style={thin, gray, -stealth}}

\usetikzlibrary{intersections}
\usetikzlibrary{angles,quotes}
\newcommand{\crosslines}[2]{
	\begin{scope}[shift={(#1,#2)}]
		\draw[-Stealth, thick] (0,0) -- (4,0);    
		\draw[Stealth-Stealth, thick] (0,3) -- (0,-3);    
			\end{scope}
}

\newcommand{\rotatedcrosslines}[3]{
	\begin{scope}[shift={(#1,#2)}]      
		\begin{scope}[rotate=#3]          
			\draw[-Stealth, thick] (0,0) -- (4,0);
			\draw[Stealth-Stealth, thick] (0,3) -- (0,-3);
		\end{scope}
	\end{scope}
}

\newcommand{\qoppa}[3]{
	\begin{scope}[shift={(#1,#2)}] 
		\begin{scope}[rotate=#3]     
			\draw[thick] (0,0) circle [radius=1cm];
			\draw[-Stealth, thick] (1,0) -- (4,0);
		\end{scope}
	\end{scope}
}

\newcommand{\anglerays}[1]{
	\pgfmathsetmacro{\cx}{cos(#1)}
	\pgfmathsetmacro{\cy}{sin(#1)}
	
	\pgfmathsetmacro{\dirA}{#1 + 135}
	\pgfmathsetmacro{\dirB}{#1 + 225}
	
	\fill (\cx,\cy) circle (1pt);
	
	\draw[thick, -{Stealth}] (\cx,\cy) -- ++(\dirA:2.3);
	\draw[thick, -{Stealth}] (\cx,\cy) -- ++(\dirB:2.3);
}

\begin{document}

\theoremstyle{plain}
\newtheorem{theorem}{Theorem}
\newtheorem{corollary}[theorem]{Corollary}
\newtheorem{lemma}[theorem]{Lemma}
\newtheorem{proposition}[theorem]{Proposition}

\theoremstyle{definition}
\newtheorem{definition}[theorem]{Definition}
\newtheorem{example}[theorem]{Example}
\newtheorem{conjecture}[theorem]{Conjecture}

\theoremstyle{remark}
\newtheorem{remark}[theorem]{Remark}

\begin{center}
{\large\bf Cutting a Pancake with an Exotic Knife} \\

\vspace*{+.1in}

David O. H. Cutler\footnote{Part of this work was carried out while D.O.H.C. was a participant in the 2025 DIMACS REU program at Rutgers University, NSF grant CCF-2447342.} \\
Department of Mathematics\\
University of Michigan\\ 
Ann Arbor, MI 48109 USA\\
Email:  \href{mailto:cutler@umich.edu}{\tt cutler@umich.edu}

\vspace*{+.1in}
 
Jonas Karlsson \\
Boden, 961 42 Sweden \\
Website: \href{https://jonka364.github.io/}{jonka364.github.io}\\
Email:  \href{mailto:jonaskarls@gmail.com}{\tt jonaskarls@gmail.com}

\vspace*{+.1in}

Neil J. A. Sloane\footnote{To whom correspondence should be addressed.} \\
The OEIS Foundation Inc.\\ 
Highland Park, NJ 08904 USA\\
Visiting Scholar, Math.\ Dept. \\
Rutgers University\\
Piscataway, NJ 08854\\
Email:  \href{mailto:njasloane@gmail.com}{\tt njasloane@gmail.com}

 \vspace*{+.1in}

\begin{abstract}
In the first chapter of their classic book {\em Concrete Mathematics}, Graham, Knuth,
and Patashnik consider the maximum number of pieces that can be obtained from a pancake
by making $n$ cuts with a knife blade that is straight, or bent into a {\sf V}, 
or bent twice into a ${\sf Z}$. We extend their work by considering knives, or ``cookie-cutters'', of even
more exotic shapes, including a $k$-armed {\sf V}, a chain of $k$ connected line segments,
long-legged versions of the letters 
${\sf A}$, ${\sf E}$, ${\sf H}$, ${\sf L}$, ${\sf M}$, ${\sf T}$, ${\sf W}$, or ${\sf X}$, a convex polygon,  
a circle,  a ${\phi}$, a figure {\sf 8}, a pentagram, a hexagram, or a lollipop (or qoppa).
We also consider ``constrained'' versions of the long-legged letters
${\sf A}$, ${\sf H}$, ${\sf L}$,  ${\sf T}$, and ${\sf X}$.
In most cases we are able to determine the maximum number
of pieces, although for the constrained ${\sf A}$ and the lollipop we can only give bounds.
\end{abstract}
\end{center}

\epigraph{
\textit{Like a long-legged fly upon the stream}\\
\textit{His mind moves upon silence.}}
{W. B. Yeats, {\em Long-Legged Fly} \cite{Yea56}}

\section{Introduction}\label{Sec1}
In this article we specify certain \mabp\  of knives  and then attempt to determine the maximum number of pieces that the plane can be divided into by making cuts from that \mabs.  Examples of \mabp\ are straight knives, {\sf V}-shaped knives, and circular knives (or cookie cutters).
For example, the members of the \mabs\ of {\sf V}-shaped knives  consist of two rays emanating from a point, with an arbitrary angle between the rays. The \mabs\ of circular knives consists of all circles of nonzero radius and arbitrary center.  
For a given \mabs\ ${\sf S}$ of knives, we wish to determine the maximum number of pieces into which the affine Euclidean plane can be decomposed by using $n \ge 0$ cuts, where each cut is chosen independently from the \mabs. We denote this number by $a_{\sf S}(n)$. 

The classic example of a \mabs\ is the set of all infinite straight lines, 
which we denote by {\sf K} (for "knife").  This is the famous 
pancake-cutting problem\footnote{Also known as the pizza-cutting problem, although we 
prefer ``pancake-cutting'', because pizzas often  have a well-defined border and for our problems
the natural setting  is the infinite affine Euclidean plane, which seems more like a pancake than a pizza.}
\cite{Bar19},  \cite[pp.\ 72--73]{Com74}, \cite[Ch.\ 20]{Gar66}, \cite[\S 1.2]{GKP94}, 
\cite[Ch.\ 1]{Mit09}, \cite[Ch.\ III]{Pol54}, 
 \cite{Ste26}, \cite[p.\ 13, Problem 44]{Yag87}, for which the answer is that the maximum number of 
pieces of an infinite pancake  that can be obtained with $n$ straight  cuts is 
\beql{Eq1}
a_{\sf K}(n) = n(n+1)/2 + 1~.
\eeq
We give a proof and an explicit construction   in \S\ref{SecPG}.
Figure \ref{FigBott} shows  solutions for $1 \le n \le 5$ cuts. 

\begin{figure}[!htb]
\begin{center}
\begin{tikzpicture} [scale =1] 

\draw[Stealth-Stealth, thick] (-1,.7) -- (1, -.3);

\def\s{3}
\draw[Stealth-Stealth, thick] (-1+\s,.7) -- (1+\s, -.3);
\draw[Stealth-Stealth, thick] (-.6+\s,1) -- (0+\s, -1);

\def\s{6}
\draw[Stealth-Stealth, thick] (-1+\s,.7) -- (1+\s, -.3);
\draw[Stealth-Stealth, thick] (-.6+\s,1) -- (0+\s, -1);
\draw[Stealth-Stealth, thick] (1+\s,.9) -- (-1+\s,-.6);

\def\s{9}
\draw[Stealth-Stealth, thick] (-1+\s,.7) -- (1+\s, -.3);
\draw[Stealth-Stealth, thick] (-.6+\s,1) -- (0+\s, -1);
\draw[Stealth-Stealth, thick] (1+\s,.9) -- (-1+\s,-.6);
\draw[Stealth-Stealth, thick] (.7+\s,1) -- (-.4+\s,-1);

\def\s{12}
\draw[Stealth-Stealth, thick] (-1+\s,.7) -- (1+\s, -.3);
\draw[Stealth-Stealth, thick] (-.6+\s,1) -- (0+\s, -1);
\draw[Stealth-Stealth, thick] (1+\s,.9) -- (-1+\s,-.6);
\draw[Stealth-Stealth, thick] (.7+\s,1) -- (-.4+\s,-1);
\draw[Stealth-Stealth, thick] (1+\s,-.1) -- (-1+\s,-.3);

\end{tikzpicture}
\end{center}
\caption{The classical pancake-cutting problem: For the \mabs\ {\sf K} of infinite straight lines, 
the plane (an infinite pancake) can be cut with 1, 2, 3, 4, or 5 cuts into a maximum of 
2, 4, 7, 11, or 16 pieces, respectively.  Arrowheads indicate infinite lines.}
 \label{FigBott}
\end{figure}

For this \mabs, the numbers $(a_{\sf K}(n))_{n \ge 0}$ form entry  \seqnum{A000124} in
the {\em On-Line Encyclopedia of Integer Sequences} (or {\em OEIS}) \cite{OEIS}.
Six-digit numbers prefixed by A will always refer to entries in the OEIS.

When we began this investigation, the OEIS already contained several other sequences of this type.
The following list includes both pre-existing  and new sequences. An asterisk~${}^\ast$ indicates a  sequence 
which was not hitherto associated with a dissection problem. 
We believe that our solutions (or, in two cases, partial solutions) to the problems indicated 
by the asterisks are new.

The sequences, associated \mabp\ ${\sf S}$, and section numbers if they are mentioned in this article, are as follows. 
The list is ordered by A-number, which is roughly the chronological order in which they were added to the OEIS
(and not the order in which they were discovered, which would be difficult to determine).

\seqnum{A000124} (${\sf S}$ is the \mabs\ {\sf K} of knife cuts, \S\S\ref{Sec1}, \ref{SecPG}),
\seqnum{A000125} (the \mabs\ of planes in three dimensions), 
\seqnum{A046127} (spheres in three dimensions), 
\seqnum{A051890} (ellipses in the plane), 
\seqnum{A058331} (hyperbolas in the plane),
\seqnum{A069894} (squares or convex quadrilaterals, \S\ref{SecSLLL}), 
\seqnum{A077588} (triangles, \S\ref{Sec1}), 
\seqnum{A077591} (concave quadrilaterals),
\seqnum{A080856}$^\ast$ (4-chains, \S\ref{SeckC}),
\seqnum{A084849}$^\ast$ (constrained $\overline{{\sf X}}$'s, \S\ref{SeccX}),
 \seqnum{A117625} (constrained ${\overline{\sf Z}}$'s or zig-zags, \S\ref{SecZMW}), 
 \seqnum{A125201}$^\ast$ (constrained ${\overline{\sf W}}$'s or ${\overline{\sf M}}$'s, \S\ref{SecZMW}), 
\seqnum{A130883} ({\sf V}'s, \S\ref{SeckV}), 
\seqnum{A140063} (constrained  $\overline{\sf H}$'s or  $\overline{\phi}$'s, \S\ref{SeccH}),
\seqnum{A140064}$^\ast$ (${\sf A}$'s, ${\sf E}$'s, 3-armed {\sf V}'s,  3-chains, Wu's (\wu),
(\begin{tikzpicture} [scale = 0.13]
\draw[thick] (0,0) -- (0,2);
\draw[thick] (-1,0) -- (-1, 2);
\draw[thick] (1,0) -- (1, 2);
\draw[thick] (-1,0) -- (1, 0);
\end{tikzpicture}), 
\S\S\ref{SeckV}, \ref{SecAG}),  
\seqnum{A143689}$^\ast$ (constrained $\overline{{\sf L}}$'s, \S\ref{SecSLLL}), 
\seqnum{A383464}$^\ast$ (4-armed {\sf V}'s, \S\ref{SeckV}), 
\seqnum{A383465}$^\ast$ (5-chains or {\sf H}'s, \S\S\ref{SeckC}, \ref{SeccH}), 
\seqnum{A383466}$^\ast$ (pentagrams, \S\ref{SecPentagram}),
\seqnum{A386477}$^\ast$ (hexagrams,  \S\ref{SecPentagram}),
\seqnum{A386480} (circles (${\sf O}$), \S\ref{SecOO}),
\seqnum{A386485}$^\ast$ (pentagons, \S\ref{SecPolygon}),
\seqnum{A386486}$^\ast$ ({\sf 8}'s, \S\ref{Sec88}),
\seqnum{A389614}$^\ast$ (constrained  $\overline{\sf T}$'s, \S\ref{SeccT}),
\seqnum{A389624}$^\ast$ (lollipops or qoppa (\qoppaletter), \S\ref{SecLP}),
\seqnum{A393442} (unions of circles and lines, \S\ref{SecSteiner}),
\seqnum{A393448}$^\ast$ (wynns (\runewynn), \S\ref{SecWynn}),
and
\seqnum{A397182}$^\ast$ (constrained  $\overline{\sf A}$'s, \S\ref{SeccA}).

The second and third entries in the above list, \seqnum{A000125} and \seqnum{A046127}, are the only three-dimensional \mabp . From this point on, we will only consider two-dimensional \mabp.

We will give references to the OEIS throughout this article, for several reasons.
One is that much of the previous work on these
problems first appeared in the OEIS rather than in a journal. We would appreciate 
hearing of any references we have overlooked, especially for the more exotic \mabp.
Another reason is that the OEIS entries include the initial terms of the
sequence,  additional formulas and  illustrations, references, etc.

As we shall see, the OEIS also led us to make two unexpected discoveries.
First, the sequences $(a_{\sf S}(n))_{n \ge 0}$ where ${\sf S}$ is any one of three \mabp, long-legged {\sf A}'s,
Wu's, or 3-chains, are identical (\S\S\ref{SeckV},  \ref{SeckC}, \ref{SecAG}).
Second,  the sequences $(a_{\sf S}(n))_{n \ge 0}$ where ${\sf S}$ is either the \mabs\ of octagons or 
of  concave quadrilaterals are  identical (\S\ref{SecOctagon}).
For the first we give a geometric explanation, whereas the second {\em may} be
simply a coincidence caused by the fact that the same formula arises in two different
but closely related contexts.

\section{Long-legged letters and other knives}\label{SecLLL}
Our interest in these exotic knives was kindled by two problems  in the first chapter of
Graham, Knuth, and Patashnik's {\em Concrete Mathematics} \cite{GKP94},
where, after discussing the classical pancake-cutting problem, they ask what would happen
if the infinite knives were bent into the shape of a {\sf V} with two very long limb: 
\begin{center}
\begin{tikzpicture} [scale = 0.50]
\draw[ -Stealth, thick] (-4,0) -- (4,.5);
\draw[ -Stealth, thick] (-4,0) -- (4, -.5);
\end{tikzpicture}
\end{center}
or bent twice so as to make a ``zig-zag'',  with two long parallel limbs
joined by a diagonal line segment:
\begin{center}
\begin{tikzpicture} [scale = 0.50]
\draw[ Stealth-, thick] ( -4,.5) -- (.5,.5);
\draw[thick] (.5,.5) -- (-.5,-.5);
\draw[-Stealth, thick] (-.5,-.5) -- (4, -.5);
\end{tikzpicture}
\end{center}
In order to generalize these unusual  knives, we observe that they may be 
regarded as ``long-legged'' versions of the upper-case letters 
{\sf V} and ${\sf Z}$, and of course an infinite straight knife can be regarded as a long-legged upper-case ${\sf I}$.
(The term was suggested by W.\ B.\ Yeats's poem in praise of silence,
 {\em Long-Legged Fly} \cite[pp.\ 327--328]{Yea56}.)

This led us to consider \mabp\ of  long-legged versions of other letters.  
Informally, a long-legged version of a letter
is obtained by replacing any protruding arms or legs by (semi-infinite) rays,
indicated by arrow-heads in the figures,
allowing arbitrary angles between these rays.
For example, there is no restriction on the angle between the two arms of a long-legged {\sf V}.
We will, of course, be careful to give a precise definition for each \mabs.  

There are interesting long-legged
versions of ${\sf A}, {\sf E}, {\sf H}, {\sf L}, 
{\sf M}, {\sf T}, {\sf W}$, and ${\sf X}$.
Curved letters 
(${\sf B}, {\sf C}, {\sf D},  {\sf G}, {\sf P}, {\sf Q}, {\sf R}, {\sf S}$, $\ldots$) do not seem appropriate 
for this treatment, since as we shall see, 
once curved knives are permitted (without any bound on their curvature),
it is easy to create arbitrarily many regions with just one or two cuts. 
We do include four curved shapes,  
$\overline{\phi}$ (\S\ref{SeccH}),
${\sf O}$ (\S\ref{SecOO}),
{\sf 8} (\S\ref{Sec88}), and the lollipop \qoppaletter\  (\S\ref{SecLP}),
that are are based on just one or two circles.
 
A knife with  the shape of a long-legged ${\sf A}$ initially looked challenging, 
although as we shall see in \S\ref{SecAG}, it has a  surprisingly simple solution.
On the other hand, the constrained long-legged 
$\overline{\sf A}$  (\S\ref{SeccA}) appears to be genuinely hard, and it and the lollipop (or qoppa) \mabs\ 
are the two \mabp\ where we do not even have a conjectured  solution (see Tables~\ref{TableAbounds}
and \ref{TableQbounds}) [July 10, 2026: It appears that the lollipop is close to being solved---see \S\ref{SecLP}.].

Long-legged versions of some letters  
(${\sf A}, {\sf H}, {\sf L}, {\sf T},  {\sf X}$, among others) 
are potentially interesting, but because our rules allow arbitrary angles 
between the extended legs, they would be indistinguishable from other letters.
An ${\sf L}$ would be the same as a {\sf V}, for example.
We therefore define \mabp\ of ``constrained'' long-legged letters, indicated by a bar over the letter,
which impose additional conditions.  Examples (which include the zig-zag knife 
$\overline{\sf Z}$ from \cite{GKP94}) will be found in Sections \ref{SecLX} and \ref{SecTA}.
 
 These constrained letters appear to be more difficult to analyze, and we will therefore 
 separate the \mabp\  we consider  into two main  {\em genera}, adopting a term from biology.
 We are working in the affine Euclidean plane, which is acted on by the
 full two-dimensional affine group $\affineGroup$, generated by linear transformations (which
 need not be isometries), followed by translations.  If ${\sf S}$ is a \mabs\ of knives, 
 we define its {\em symmetry group} 
 $\sG({\sf S})$ to be the largest subgroup of $\affineGroup$ which preserves ${\sf S}$, 
 that is, has the property that, for any $g \in \sG({\sf S})$ and any $s \in {\sf S}$,  we have $g(s) \in {\sf S}$.
 
 The two genera are  the {\em affine genus}, consisting of those \mabp\ ${\sf S}$
 for which $\sG({\sf S}) = \affineGroup$,  and the {\em similarity genus}, consisting
 of those \mabp\ for which $\sG({\sf S}) = \simGroup$, the group of all similarity transformations, generated
 by  translations, reflections, rotations, and dilations (or rescaling).  We further split each of these genera into two subgenera, depending on whether or not the
 symmetry group $\sG({\sf S})$ acts transitively on ${\sf S}$. Table~\ref{TableG} shows
 how our principal  \mabp\ fall into these four subgenera.
 (These \mabp\ will be defined in the following sections.)  We hope that this table, combined with the list
 given in \S\ref{Sec1}, will bring order into what might seem like a bewildering array of shapes.
\begin{table}[!htb]
$$
\begin{array}{c|c|c}
\rm{genus}             & \affineGroup                          & \simGroup  \\
\hline  
\rule{0pt}{14pt}
{\rm transitive}   &  {\sf K}\, (\ref{SecPG}),    \begin{tikzpicture} [scale = 0.25] \draw[fill, ultra thick] (0,0) circle (.1); \draw[thick] (0,0) -- (0, -1); \end{tikzpicture} \, (\ref{SecH}),  
 {\sf A}\,  (\ref{SecAG}), 
{\runewynn} \, (\ref{SecWynn})&     
\overline{\sf L},    \overline{\sf X},   \overline{\sf H},  \overline{\phi}\, (\ref{SecLX}),
 \overline{\sf T}\, (\ref{SecTA}),
   {\sf O},    {\sf 8},  {\qoppaletter}\, (\ref{SecFS2}) \\
{\rm intransitive} & k{\sf V}\,  (\ref{SeckV}),  k{\sf C}\, (\ref{SeckC}) &
\overline{\sf Z},   \overline{\sf W}\, (\ref{SecZMW}), \overline{\sf A}\, (\ref{SecTA}) \\
\end{array}             
$$                      
  \caption{ Principal \mabp\ of knives classified into four subgenera,
  together with their section numbers.}\label{TableG}
\end{table}

Of course there are other important subgroups of $\affineGroup$ that we could have focused on. If our knives
were circles of a fixed radius, the appropriate subgroup 
would be the group $\isoGroup$ of isometries of the plane, which is strictly contained in the group of similarities: 
$\isoGroup \subset \simGroup \subset \affineGroup$.
Or we might wish to study only orientation-preserving transformations.
But neither of these groups will be considered here.

As their symbols suggest,  the \mabp\  ${\sf O}$, ${\sf 8}$,  $\overline{{\phi}}$, and  \qoppaletter,  are defined in terms of circles, which puts them in the similarity genus, since affine maps do not take circles to circles. 
Likewise, as their symbols suggest, the \mabp\ $\overline{\sf L}$,  $\overline{{\sf X}}$, and $\overline{\sf T}$ 
involve perpendicular lines,  and so are also not affine.

For a transitive \mabs\  of knives, any single example will serve as a {\em holotype}, or typical member.
We will call such a \mabs\ {\em holotypic}, and it can be specified simply by giving its symmetry group and  the holotype.

Since $\affineGroup$ is 3-transitive, the family of non-degenerate triangular knives is a holotypic
member of the transitive affine subgenus. In contrast, the family of {\em right-angled}  triangular knives 
is an intransitive, or {\em heterotypic}, member of the similarity subgenus.

For the family of non-degenerate quadrilaterals (or those using any polygon
with more than three vertices),  the symmetry group is $\affineGroup$, which now  acts intransitively, 
since there are not enough transformations in $\affineGroup$
to take one quadrilateral to any other.
Similarly,  for $k>2$, $k$-chains and $k$-armed {\sf V}'s also have too many degrees of freedom to be holotypic. In a heterotypic \mabs\ ${\sf S}$, there would be a holotype for each orbit of the action of  $\sG({\sf S})$.

Generally speaking, the bigger the group $\sG({\sf S})$, the easier
it is to find $a_{\sf S}(n)$,  because there are more possibilities for the  next knife. 
Intransitive \mabp\  also give us more freedom. 

\section{Further remarks on the families }\label{SecFurther}
For a given \mabs\ ${\sf S}$ of knives, our strategy for attempting to determine  $a_{{\sf S}}(n)$
is to select $n$  members of ${\sf S}$ and analyze 
the  planar graph  $\Gamma_{\sf S}(n)$ formed by their union. 
We describe this graph in \S\ref{SecPG},  using the classical case when the knives are lines as an illustration.

Subsequent sections will then discuss other shapes, roughly in order of increasing complexity,
beginning in \S\ref{SecH} with the half-line or ray (\S\ref{SecH}).
To emphasize the difference between a line and a ray, we refer to a ray as a {\em hatpin},
with symbol \begin{tikzpicture} [scale = 0.25] \draw[fill, ultra thick] (0,0) circle (.1); \draw[thick] (0,0) -- (0, -1); \end{tikzpicture}
(see Fig.\ \ref{FigHatpinWu} below for some drawings  of hatpins).

The next-simplest  shape after the hatpin is the {\sf V},
and in~\S\ref{SeckV}  we will analyze a long-legged $k$-armed {\sf V} for $k \ge 2$.
The three-armed {\sf V}, \wu,  is particularly interesting. It is called a Wu, 
and is included in the Unified Canadian Aboriginal Syllabics section of Unicode, 
where it has the  symbol U+15D0. The Wu has also been called a Boolean OR with a middle stem (U+2A5B).
The numbers of regions obtained from a long-legged ${\sf E}$ 
or  Shah 
(\begin{tikzpicture} [scale = 0.13]
\draw[thick] (0,0) -- (0,2);
\draw[thick] (-1,0) -- (-1, 2);
\draw[thick] (1,0) -- (1, 2);
\draw[thick] (-1,0) -- (1, 0);
\end{tikzpicture})
turn out to be the same as for a Wu, so we simply define  long-legged ${\sf E}$'s and 
Shah's  to be long-legged Wu's.

The sections following the multi-armed ${\sf V}$ 
 are devoted to the  multi-armed $k$-chain (\S\ref{SeckC}), the  
long-legged  ${\sf A}$ (\S\ref{SecAG}), 
and the wynn,   a Middle English rune  (\S\ref{SecWynn}).
These are followed by the constrained long-legged letters 
$\overline{\sf Z}$, $\overline{\sf W}$, and $\overline{\sf M}$ (\S\ref{SecZMW}),
and  in \S\ref{SecLX} by  $\overline{{\sf L}}$,
$\overline{{\sf X}}$, $\overline{\sf H}$, and  $\overline{{\phi}}$.
Our solution to the $\overline{\sf H}$ \mabs\ given in \eqref{EqcH} led us to
D. Kinsella's work on the $\overline{{\phi}}$ \mabs, where the same formula
as for $\overline{\sf H}$  appears, although without proof.
However, this formula, \eqref{EqcH}, had already been given 200 years earlier by J. Steiner \cite{Ste26},
as we briefly discuss in \S\ref{SecSteiner}.

 Two further (and more difficult) constrained letters, $\overline{{\sf T}}$ and $\overline{{\sf A}}$,  follow in \S\ref{SecTA}.
 
We study several finite shapes in 
\S\ref{SecFS}, namely  a convex $k$-sided 
polygon (\S\ref{SecPolygon}, \S\ref{SecOctagon}),
a regular pentagram and hexagram (\S\ref{SecPentagram}),
and then  three curved \mabp\, 
the circle (\S\ref{SecOO}),  figure {\sf 8}  (\S\ref{Sec88}), and ending with
the  lollipop or qoppa (\qoppaletter) (\S\ref{SecLP}).

Our methods can be applied to many shapes besides these, and could  be used to provide proofs 
for many formulas stated without proof in the OEIS entries mentioned in our initial list.

The final section lists some open problems. One in particular is worth emphasizing:
is there any connection between the present work and the classical
subject of ``Geometrical Configurations'' \cite{Gru09, HCV52, Lev29}?
Many of our graphs are in fact geometrical configurations, so it is natural to
wonder if any results from the classical literature can
be  applied to our problem.

\vspace*{+.1in}

  \begin{figure}[!htb]
 \centerline{\includegraphics[angle=0, width=3.0in]{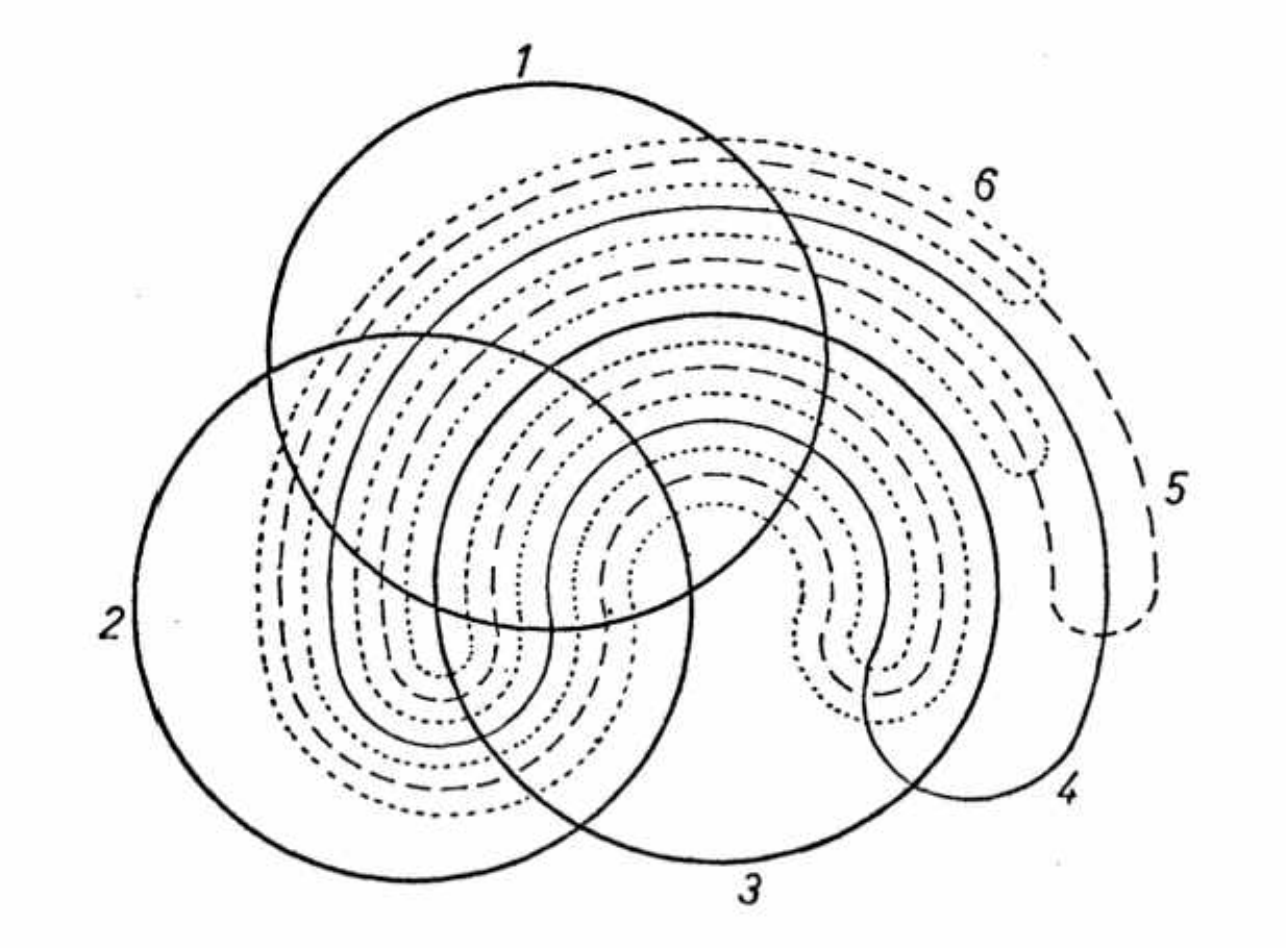}}
 \caption{A $6$-set Venn diagram (with $64$ regions) constructed 
 from six copies of a sausage shape $S$, from  Gr\"{u}nbaum \cite{Gru75}.
 (A circle is a very simple sausage.) 
 The construction easily generalizes to show that at least $2^n$ regions can be achieved 
 using $n$ copies of  a simple Jordan curve.} \label{FigVenn}
 \end{figure}

For most of the \mabp\ that we consider (and also  for most of those on the  
list at the start of this Introduction),  
 $a_{\sf S}(n)$ grows quadratically with $n$.  For triangles (\seqnum{A077588}), for example,  $a_{\sf S}(n) = 3n^2-3n+2$ for $n \ge 1$. 
It is possible, however, for $a_{\sf S}(n)$ to grow much more rapidly than this.
For example, Gr\"{u}nbaum's $n$-set Venn diagram construction \cite{Gru75, RuWe05}  
(see Fig.\ \ref{FigVenn} for $n=6$), shows that if $S$ is a simple Jordan curve,  
a snake or sausage shape, for example, we can achieve $a_{JC}(n) \ge 2^n$. 
In fact, $a_{JC}(2)$ is already unbounded, since it can be made as large as we please
by intertwining two snakes  (as  in  drawings of a caduceus).
Of course the  very definition of a Jordan curve implies  $a_{JC}(1) = 2$.
But even $a_{\sf S}(1)$ itself can be arbitrarily large if $S$ is a twisted sausage (Fig.\ \ref{FigSaus}).
To avoid such complications, we will mostly restrict ourselves to polygonal shapes. 
\begin{figure}[!ht]
        \begin{center}
                \begin{tikzpicture}[rotate=0]  [scale = 2]
\draw[thick] (0,0) -- (9,0);
\draw[thick, rounded corners] (0,0) to [out=180, in=180] (0.5,-0.325) 
   to [out=0, in=270] (1,0)
   to [out=90, in=90] (2,0)
   to [out=270, in=270] (3,0) 
   to [out=90, in=90] (4,0)
   to [out=270, in=270] (5,0) 
   to [out=90, in=90] (6,0)
   to [out=270, in=270] (7,0) 
   to [out=90, in=90] (8,0)
   to [out=270, in=180] (8.5,-0.325) 
   to [out=0, in=0] (9,0);
                \end{tikzpicture}
        \end{center}
\caption{A single twisted sausage shape can produce arbitrarily many regions.}
\label{FigSaus}
\end{figure}

\section{Notation}\label{SecNotn} 
Many of our knife \mabp\ are named after letters of the alphabet, so
in order to distinguish between the names of the  \mabp\ and
the mathematical symbols used for their properties,  we will use a {\em sans serif}  font
({\sf A},  $\overline{{\sf L}}$, {\sf V}, etc.) for the \mabp.
Other \mabs\ names are  {\sf K} for the infinite knife-cuts
in the basic pancake-cutting problem,  
\begin{tikzpicture} [scale = 0.25] \draw[fill, ultra thick] (0,0) circle (.1); \draw[thick] (0,0) -- (0, -1); \end{tikzpicture}
for the hatpin family,  
$k{\sf V}$ for the $k$-armed {\sf V}'s,
$k{\sf C}$ for the $k$-chains,
$k{\sf P}$ for the $k$-sided polygons,
${\sf O}$ for  circles (it should really be $\overline{\sf O}$, but we omit the bar for aesthetic reasons),
and {\sf St} (honoring Steiner) for  knives formed from $m$ circles and $n$ lines.
An overbar ($\overline{{\sf L}}$, $\overline{{\sf X}}$, etc.) indicates a constrained shape.

For a \mabs\ ${\sf S}$, $s \in {\sf S}$ is a knife in the \mabs, 
$\Gamma_{\sf S}(n)$ is any graph formed by drawing $n$ knives $s \in {\sf S}$, and
$a_S(n)$ is the maximum number of regions in any $\Gamma_{\sf S}(n)$. 
We will say that $\Gamma_{\sf S}(n)$ is {\em optimal} if it has $a_S(n)$ regions.
Red and black dots in $\Gamma_{\sf S}(n)$ usually  indicate base and crossing nodes, respectively (these terms will be defined in the next section).

The parameters of $\Gamma_{\sf S}(n)$ will be denoted by  $V_{\sB}$, $V_{\sC}$, $V = V_{\sB} + V_{\sC}$ (numbers of base nodes, crossing nodes,  total vertices), $\Ei$, $E_f$, $E = \Ei + E_f$ (infinite edges, finite edges, total edges),
and $\Ri$, $R_f$, $R = \Ri+R_f$ (infinite regions, finite regions, total regions).
Also, $\alpha$ denotes the total number of arms in the graph,  and $\cpa$ is the 
average number of {\em crossing nodes per arm}.

For individual knives, $\sigma(S)$ is the maximum number
of self-crossings in a single $s \in {\sf S}$, 
$\kappa(S)$ is the maximum number of intersections between
$s \ne s^{\prime}  \in {\sf S}$, and
$\xi(S)$ (which will not appear until \S\ref{SeccA}) is a local intersection number. 

The exotic symbols used in this article  are:
\wu\ (Unified Canadian Aboriginal symbol Wu, \S\ref{SecFurther}),
\begin{tikzpicture} [scale = 0.25] \draw[fill, ultra thick] (0,0) circle (.1); \draw[thick] (0,0) -- (0, -1); \end{tikzpicture}
(hatpin, \S\ref{SecH}),
\runewynn\ (wynn, a Middle English rune,  \S\ref{SecWynn}), and
\qoppaletter\ (lollipop or archaic Greek letter qoppa, \S\ref{SecLP}).

\section{The planar graph and the maximum number of regions}\label{SecPG}
Given a \mabs\ ${\sf S}$, we start by thinking of the knives $s \in {\sf S}$  as planar graphs with certain vertices and edges.
We call the  vertices of $s$  {\em base nodes} and its edges {\em arms}.
A long-legged {\sf V}, for example, has a single base node and two arms. 
A $k$-chain (\S\ref{SeckC}) has $k-1$ base nodes and $k$ arms.

The planar graph $\Gamma_{\sf S}(n)$ for a drawing formed from $n$ knives $s \in S$ is then the union of 
the graphs for the individual knives, 
together with a new vertex for any point where a pair of arms cross,
and new edges that replace the portions of the arms between the new vertices.
We will refer to these new vertices as {\em crossing nodes}.
The arms of an $s \in S$ may themselves cross, and these {\em self-crossings}
are  also considered to be crossing nodes.
In some drawings the base nodes are indicated by small red circles, and the crossing nodes by small black circles.
(We have not done this in every drawing, to avoid clutter.)

We denote the set of base nodes for all $n$ knives  by $\sB$, and the degree of $v \in \sB$ by $d_v$,
or by $d_{\sB}$ if the degree is the same for all $v \in \sB$.

Since our goal is always to maximize the number of regions, we will assume
that no crossings ever involve more than two arms,
and that  no arm from one $s$ intersects a base node from a different $s$
(for otherwise a small perturbation would increase the number of regions).
Since we are constructing these graphs ourselves, it is clear that these assumptions are justified.
A similar perturbation argument is briefly mentioned in \cite[p.\ 128]{Mat02}.

Let $V_{\sB}$ denote the total number of base nodes over all $n$ knives,
$V_{\sC}$ the total number of crossing nodes,
and $V=V_{\sB}+V_{\sC}$ the total number of vertices of $\Gamma_{\sf S}(n)$. 
Let $\Ei$ and $E_f$ respectively  denote the number of infinite and finite edges, and 
$E = \Ei + E_f$ the total number of edges, and similarly let
$\Ri$, $R_f$, and $R=\Ri+R_f$  denote the numbers of 
infinite, finite, and total  regions.

If  the knives $s \in S$ are bounded and $\Ei = 0$,  Euler's formula 
\cite[p.\ 22]{Bol98}, \cite[Ch.~8]{Gru67}, 
\cite[\S2.2]{Gru72}, \cite{Lak76}  for  $\Gamma_{\sf S}(n)$ 
tells us that
\beql{Euler2}
R-E+V= 2 ~.
\eeq
This version of Euler's formula will be relevant in~\S\ref{SecFS}.
But most of the knives we consider are infinite, with $\Ei = \Ri >  0$,  in which case 
Euler's formula for the Euclidean (or affine) plane tells us that $\Gamma_{\sf S}(n)$ satisfies 
\beql{Euler1}
R-E+V= 1 ~.
\eeq
We were surprised not to find \eqref{Euler1}  in Lakatos \cite{Lak76}, 
a  book devoted to  proofs of different versions of Euler's formula,
nor in  the standard books  on finite graph theory.  But \eqref{Euler1}  is simple
to prove.  Draw a rectangle around $\Gamma_{\sf S}(n)$, large enough so that only the $\Ei$ edges
cross it. The rectangle itself (the frame) therefore has $4+\Ei$ edges and $4+\Ei$ vertices. 
Discard all portions of the $\Ei$ edges that are outside the frame. The resulting finite graph
has $E + 4 + \Ei$ edges, $V+4+ \Ei$ vertices,
and by~\eqref{Euler2} has $(E+4+\Ei) - (V+4+\Ei) + 2 = E-V+2$ regions.
But this includes the exterior of the rectangle, so inside the rectangle there are $V-E+1$ regions,
which establishes \eqref{Euler1}.

Let us assume then that the knives in ${\sf S}$  are connected graphs formed from straight line segments, some of which may be infinite.
What greatly simplifies  the analysis of $\Gamma_{\sf S}(n)$
is that there is a second equation  linking 
$E$ and $V$, which we obtain by counting edge-vertex incidences in two ways.
On the left side of this equation we sum the degrees~$d_v$ of all the base 
nodes $v \in \sB$
and add four times the number of crossing nodes. This will count each edge twice,
except for the $\Ei$ infinite edges, which will only be counted once, and so we add 
an extra term $\Ei$ to the left side. On the right side we simply count each edge twice.
The edge-vertex count  then tells us that
\beql{EV0}
\sum_{v \in \sB} d_v + 4 V_{\sC} + \Ei = 2(\Ei + E_f)~,
\eeq
that is,
\beql{EqEV}
E_f = 2 V_{\sC} + \frac{1}{2} \sum_{v \in \sB} d_v  - \frac{1}{2} \Ei~.
\eeq
By combining this with Euler's formula \eqref{Euler1}, 
we obtain our basic equation, which is
\beql{EqEV2}
R = V_{\sC} + \frac{1}{2} \sum_{v \in \sB} d_v   - V_{\sB} + \frac{1}{2}\Ei +1~.
\eeq
 
 On the right side, only $V_{\sC}$ depends on how the knives are placed in the plane, the other 
 quantities being determined by ${\sf S}$ and $n$, and so we deduce that, for the \mabp\ we are considering,
\begin{align}\label{EqL3}
 & \mbox{to~maximize~the~number~of~regions,~it~is~necessary} \nonumber \\
 & \mbox{and~sufficient~to~maximize~the~number~of~crossings.}
\end{align}

Note that this conclusion holds both for finite and infinite knives: the constant term in the formula is different in the two cases, but its value is irrelevant.

We now illustrate this by studying  the classic  pancake-cutting problem mentioned in the Introduction,
for which the knives $s \in {\sf K}$  are infinite straight lines, such as  the line $y = 0$ (the ``holotype''),  
or any of  its images 
under the affine group $\affineGroup$. 
The knives $s$ have one arm each and  no  base vertices. 
We form the pancake graph $\Gamma_{\sf K}(n)$ by drawing $n \ge 2$ such lines.
 We may assume that each line crosses at least one other line, so $V_{\sB} = 0$
 and $\Ei = 2n$.  Equation \eqref{EqEV2} then says that
\beql{EqL2}
R = V_{\sC} + n + 1~.
\eeq
We maximize $V_{\sC}$  by making every pair of lines intersect, giving
$V_{\sC} = \binom{n}{2}$, and so the maximum number of regions is
$R = a_{\sf K}(n) = \binom{n}{2}+n+1$, establishing \eqref{Eq1}.  

The underlying reason that \eqref{EqL3} holds  is that, because our shapes 
are built out of straight lines, locally the graphs have a tendency  to look 
like a square mesh, with $E  \approx 2V$
(cf. \eqref{EqEV}, \eqref{EqkV0}). This also explains why the 
associated sequences
$(a_{\sf S}(n))$ have only quadratic growth.

\begin{figure}[!htb]
        \begin{center}
                \begin{tikzpicture}[rotate=0]  [scale=1]
\newcommand*{\TickSize}{2pt}%

\begin{scope}[scale=0.7]
\draw [axis line style] (-1,0) -- (9,0);
\draw [axis line style] (0,-3) -- (0,5);

\foreach \x in {0,1,...,9} {%
    \draw ($(\x,0) + (0,-\TickSize)$) -- ($(\x,0) + (0,\TickSize)$)
        node [below] {$\x$};
}

\foreach \y in {-3,-2,...,-1} {%
    \draw ($(0,\y) + (-\TickSize,0)$) -- ($(0,\y) + (\TickSize,0)$)
        node [left] {$\y$};
}
\foreach \y in {1,2,...,5} {%
    \draw ($(0,\y) + (-\TickSize,0)$) -- ($(0,\y) + (\TickSize,0)$)
        node [left] {$\y$};
}
\draw[Stealth-Stealth, thick] (-1,0) -- (9,0);
\draw[Stealth-Stealth, thick] (-1,-.45) -- (9,1.95);
\draw[Stealth-Stealth, thick] (-1,-1.67) -- (9,3.90);
\draw[Stealth-Stealth, thick] (0,-2.80) -- (8,4.90);
\draw[Stealth-Stealth, thick] (2,-3.15) -- (7,5.00);
\draw[Stealth-Stealth, thick] (4.15,-3.00) -- (6.4,5.00);
\draw[Stealth-Stealth, thick] (6,-3.00) -- (6,5.00);
\end{scope}
                \end{tikzpicture}
        \end{center}
 \caption{An optimal pancake graph $G_{\sf K}(7)$   defined by \eqref{EqPancake}
 with $n=6$ and $\theta = 15^{\circ}$. 
 The $t=0$ line is the $x$-axis, and the $t=6$ line is $x=6$.}
 \label{FigKnife2}
\end{figure}

Since the optimal pancake graph $\Gamma_{\sf K}(n)$ is the basis for other constructions, 
we summarize its  properties: $V = V_{\sC} = \binom{n}{2}$,
$\Ei = \Ri = 2n$, $E_f = n(n-2), E = n^2$,
$R_f = \binom{n-1}{2}$, and $R = a_K(n) = (n^2+n+2)/2$ (\seqnum{A000124}).
Our very first figure, Figure~\ref{FigBott} above,  shows some {\em ad hoc} drawings of 
these graphs. The following is an explicit construction of an optimal pancake graph 
with $N+1 \ge 2$ lines. Let $\theta = \pi/(2N)$, and take the lines
to have equations 
\beql{EqPancake}
(x-t) \sin t \theta =  y \cos t \theta \quad \mbox{~for~} t = 0, 1, \ldots N~.
\eeq
Figure~\ref{FigKnife2} illustrates the construction  of a $\Gamma_{\sf K}(7)$, taking $N=6$ .

For use in \S\ref{SecFS} and \S\ref{SecFS2} we give the analogue of our basic equation \eqref{EqEV2} that 
applies if the knives $s \in {\sf S}$ are finite graphs.  Since  $\Ei=0$, the equation is now 
\beql{EqEVf}
R = V_{\sC} + \frac{1}{2} \sum_{v \in \sB} d_v   - V_{\sB} + 2~.
\eeq

\vspace*{+.1in}

We end this section by defining three parameters associated with a \mabs\  ${\sf S}$ that are helpful
when  searching for optimal graphs. First, if the definition of ${\sf S}$
permits, it may be possible to move the arms of ${\sf S}$ so as to increase the number of self-crossings (this applies to the $k$-chain, for example). 
We let $\sigma({\sf S})$ denote the maximum number
of self-crossings in a single $s \in {\sf S}$.

Second, $\kappa({\sf S})$ will denote the maximum number of intersections between
two knives \mbox{$s \ne s^{\prime}  \in {\sf S}$}. With $n$ copies of ${\sf S}$, an upper bound on the total
number of crossings is therefore 
\beql{EqVx}
V_{\sC} \le n \,\sigma({\sf S}) + \binom{n}{2} \kappa({\sf S})~.
\eeq
If equality holds in \eqref{EqVx},  the graph is optimal.

In many of the examples studied in this article, we are able to achieve equality in \eqref{EqVx}.
The chief difficulty, once we have determined $\kappa({\sf S})$, lies in placing the $n$ knives so that every pair intersect in $\kappa({\sf S})$ points. We run into difficulties if,
once one $s$ has been positioned, there is only a small range of
possibilities for the second (and therefore  {\em all}) subsequent copies.
This is the case for the constrained letters $\overline{\sf T}$ (\S\ref{SeccT}) and $\overline{\sf A}$ (\S\ref{SeccA}),
the figure ${\sf 8}$ (\S\ref{Sec88}), and the lollipop (\S\ref{SecLP}).

Third, by assumption, each crossing node is at the intersection of exactly two  
arms, and so the average number of crossings per arm  (abbreviated $\cpa$) is
\beql{Eqcpa1}
\cpa = \frac{2 V_{\sC}}{\alpha(n)} \le \frac{n \,\sigma({\sf S}) + \binom{n}{2} \kappa({\sf S})}{\alpha(n)}~,
\eeq
where $\alpha(n)$ is the total number of arms.
The $\cpa$ is not always an integer, even for optimal graphs 
(the $k$-chains for $k \ge 3$, for example). But if the right side of \eqref{Eqcpa1}
{\em is} an integer, the $\cpa$  can be useful when searching for optimal graphs.
If in the graph we are currently trying to optimize there is an 
arm with fewer than $\cpa$  crossings, then either that arm needs to be changed, 
or  it  is a signal that  the bound in \eqref{Eqcpa1} cannot be achieved.

\section{The hatpin}\label{SecH}
We define a {\em hatpin}
to be a planar graph consisting of a distinguished point, 
its head, with a ray emanating from it.  The \mabs\ of all hatpins will be denoted by 
\begin{tikzpicture} [scale = 0.25] \draw[fill, ultra thick] (0,0) circle (.1); \draw[thick] (0,0) -- (0, -1); \end{tikzpicture}.   
Its group is $\affineGroup$, which acts transitively, and for holotype we may take the ray
 $ \{x=0,\, y \le 0\}$ with head $ (0,0)$. 
A graph 
$\Gamma_{\begin{tikzpicture} [scale = 0.25] \draw[fill, ultra thick] (0,0) circle (.1); \draw[thick] (0,0) -- (0, -1); \end{tikzpicture}}(n)$ 
formed by drawing $n$ hatpins has $V_{\sB} = n$ base nodes, 
$d_B=1$,  and $\Ei=n$, so \eqref{EqEV2} gives $R = V_{\sC} + 1$, which we maximize by making all the hatpins cross each other, giving  
$a_{\begin{tikzpicture} [scale = 0.25] \draw[fill, ultra thick] (0,0) circle (.1); \draw[thick] (0,0) -- (0, -1); \end{tikzpicture}}(n) = \binom{n}{2} + 1$.
The left-hand figure in Figure~\ref{FigHatpinWu} below shows a 
$\Gamma_{\begin{tikzpicture} [scale = 0.25] \draw[fill, ultra thick] (0,0) circle (.1); \draw[thick] (0,0) -- (0, -1); \end{tikzpicture}}(4)$.
In some drawings we have indicated the heads of the hatpins by small red circles.

Since the optimal hatpin graph $\Gamma_{\begin{tikzpicture} [scale = 0.25] \draw[fill, ultra thick] (0,0) circle (.1); \draw[thick] (0,0) -- (0, -1); \end{tikzpicture}}(n)$ is also used in other constructions, 
we note that it has parameters 
$V_{\sB} = n,  d_B = 1, V_{\sC} = \binom{n}{2}, V = \binom{n+1}{2}$,
$\Ei = \Ri = n$, $E_f = n(n-1), E = n^2$,
$R_f = \binom{n-1}{2}$, and $R = a_{\sf K}(n-1) = (n^2-n+2)/2$.

\begin{figure}[!htb]
        \begin{center}
                \begin{tikzpicture}[rotate=0]  [scale=1]

\begin{scope}[rotate=-7] [scale=0.5]
\draw[ -Stealth, thick, black] (1.4,4.5) -- (0.6,0.4);
\draw[ -Stealth, thick, black] (1.4,4.5) -- (1.8,0.4);
\draw[ -Stealth, thick, black] (1.4,4.5) -- (3.0,0.4);
\draw[ -Stealth, thick, blue] (-0.3,2.2) -- (5,3.9);
\draw[ -Stealth, thick, blue] (-0.3,2.2) -- (5,2.8);
\draw[ -Stealth, thick, blue] (-0.3,2.2) -- (5,1.7);
\draw[red, ultra thick, fill] (1.4,4.5) circle [radius=.08];
\draw[red, ultra thick, fill] (-.3,2.2) circle [radius=.08];
\end{scope}
                \end{tikzpicture}
        \end{center}
\caption{A graph $G_{3{\sf V}}(2)$ with 14 regions, formed from two 3-armed {\sf V}'s (or Wu's).}
 \label{FigWu2}
\end{figure}

 \begin{table}[!htb]
$$
\begin{array}{c|rrrrrrrrrr}
k \backslash n & 0 & 1 & 2 & 3 & 4 & 5 & 6 & 7 & 8 & ... \\
\hline
0 &   1 &  1 &  1 &  1 &  1 &  1 &  1 &  1 &  1 &  ... \\
1 &   1 &  1 &  2 &  4 &  7 &  11 &  16 &  22 &  29 &  ... \\
2 &   1 &  2 &  7 &  16 &  29 &  46 &  67 &  92 &  121 &  ... \\
3 &   1 &  3 &  14 &  34 &  63 &  101 &  148 &  204 &  269 &  ... \\
4 &   1 &  4 &  23 &  58 &  109 &  176 &  259 &  358 &  473 &  ... \\
5 &   1 &  5 &  34 &  88 &  167 &  271 &  400 &  554 &  733 &  ... \\
6 &   1 &  6 &  47 &  124 &  237 &  386 &  571 &  792 &  1049 &  ... \\
7 &   1 &  7 &  62 &  166 &  319 &  521 &  772 &  1072 &  1421 &  ... \\
\end{array}
$$
  \caption{Table of $a_{k{\sf V}}(n)$, the maximum number of regions formed in the plane by drawing $n$ long-legged $k$-armed {\sf V}'s.  The entries are given by \eqref{EqkV1}. Rows 1-4 are 
  \seqnum{A000124}, \seqnum{A130883}, \seqnum{A140064}, 
  \seqnum{A383464},
  column 2 is \seqnum{A008865}, the main diagonal is \seqnum{A393441},
  and the OEIS entry for the array itself is \seqnum{A386481}.}\label{Table1}
\end{table}

\section{The \texorpdfstring{$k$}{k}-armed {\sf V}}\label{SeckV}
Our next shape is a  long-legged $k$-armed {\sf V},
which we define to consist of a distinguished point, its tip,  
with $d_B = k \ge 1$ rays, or arms,  emanating from it,
with arbitrary angles between the arms. We denote the \mabs\ of all such knives
by $k{\sf V}$.  Its group is $\affineGroup$, but if $k>2$ it does not act transitively:
this is a heterotypic \mabs.
A  graph 
$\Gamma_{k{\sf V}}(n)$ formed by drawing $n \ge 0$ $k$-armed ${\sf V}$'s  has 
$V_{\sB} = n$ base nodes and  $nk$ arms,
with $\Ei = \Ri = nk$. 
For maximizing the number of regions, it turns out to be best to take the arms to be distinct, 
and  forming a fan with a small angular spread, as in Figs.\ \ref{FigWu2} and \ref{FigHatpinWu}.

From Equations~\eqref{EqEV}  and \eqref{EqEV2} we obtain
\beql{EqkV0}
E_f = 2 V_{\sC}, ~~ R = V_{\sC} + n(k-1) + 1~.
\eeq
The crossing number $\kappa(k{\sf V})$ is clearly $k^2$.
We maximize $V_{\sC}$ by making every pair of $k{\sf V}$'s cross each other in $k^2$ points, 
for a total value of
$V_{\sC} =  \binom{n}{2} k^2$,
leading to the maximum of 
\beql{EqkV1}
R = a_{k{\sf V}}(n) =  \binom{n}{2} k^2  +n(k-1) +1\quad  (k \ge 1, n \ge 0)
\eeq
regions.\footnote{The $k=3$ case of \eqref{EqkV1}, $(9n^2-5n+2)/2$ (\seqnum{A140064}),  was first discovered by Edward Xiong, Jonathan Pei, and David O. H. Cutler on June 24 2025 (see Acknowledgments).}  This is easily achieved by taking 
a $\Gamma_{\begin{tikzpicture} [scale = 0.25] \draw[fill, ultra thick] (0,0) circle (.1); \draw[thick] (0,0) -- (0, -1); \end{tikzpicture}}(n)$ hatpin graph  
and replacing the hatpins by  long narrow $k$-armed ${\sf V}$'s.
(see Fig.\ \ref{FigHatpinWu}, which illustrates the construction for four 3-armed ${\sf V}$'s).
We note that for $k=2$ our counts match those of \cite{GKP94}, as they should. Those authors give a different explicit construction. They also note (as their Exercise 19) that imposing a minimal angle
between the arms of the ${\sf V}$  makes this maximal count unattainable; whichever construction we use,  the ${\sf V}$'s must necessarily get thinner and thinner.  Table~\ref{Table1} shows the values of $a_{k{\sf V}}(n)$ for
$0 \le k \le 7$ and $0 \le n \le 8$.

\begin{figure}[!htb]
        \begin{center}
                \begin{tikzpicture}[rotate=0]  [scale=1]

\begin{scope}[rotate=0] [scale=1]
\def\s{4}
\draw[ Stealth-,  thick, black] (0.6,1.5) -- (3.85,2.9);
\draw[ Stealth-,  thick, black] (0.5,1.85) -- (4.1,1.5);
\draw[ -Stealth,  thick, black] (0.4,4.3) -- (3.7,.8);
\draw[ -Stealth,  thick, black] (1.9,5.0) -- (1.9,0.1);
\draw[red, ultra thick, fill] (.4,4.3) circle [radius=.08];
\draw[red, ultra thick, fill] (1.9,5) circle [radius=.08];
\draw[red, ultra thick, fill] (3.85,2.9) circle [radius=.08];
\draw[red, ultra thick, fill] (4.1,1.5) circle [radius=.08];
\draw[ thick, black] (\s+0.6,1.5) -- (\s+3.85,2.9);
\draw[ Stealth-,  thick, black] (\s+0.6,1.4) -- (\s+3.85,2.9);
\draw[ thick, black] (\s+0.6,1.3) -- (\s+3.85,2.9);
\draw[ thick, black] (\s+0.5,1.85) -- (\s+4.1,1.5);
\draw[ Stealth-,  thick, black] (\s+0.5,1.95) -- (\s+4.1,1.5);
\draw[ thick, black] (\s+0.5,2.05) -- (\s+4.1,1.5);
\draw[ thick, black] (\s+0.4,4.3) -- (\s+3.7,.8);
\draw[ -Stealth, thick, black] (\s+0.4,4.3) -- (\s+3.7,1.0);
\draw[ thick, black] (\s+0.4,4.3) -- (\s+3.7,1.2);
\draw[ thick, black] (\s+1.9,5.0) -- (\s+1.9,0.1);
\draw[ -Stealth, thick, black] (\s+1.9,5.0) -- (\s+2.05,0.1);
\draw[ thick, black] (\s+1.9,5.0) -- (\s+2.2,0.1);
\draw[red, ultra thick, fill] (\s+.4,4.3) circle [radius=.08];
\draw[red, ultra thick, fill] (\s+1.9,5) circle [radius=.08];
\draw[red, ultra thick, fill] (\s+3.85,2.9) circle [radius=.08];
\draw[red, ultra thick, fill] (\s+4.1,1.5) circle [radius=.08];
\node[right, blue] at (2.1,5.0) {I};
\node[right, blue] at (2.0,3.8) {II};
\node[right, blue] at (2.5,2.2) {III};
\node[right, blue] at (2.30,1.2) {IV};
\node[left, blue] at (1.75,0.8) {V};
\draw[->, thick,blue] (2.44,1.4) -- (2.20,2.00);

\end{scope}
                \end{tikzpicture}
        \end{center}
\caption{Transforming a hatpin graph $G_{\sf H}(4)$ to an optimal 
$G_{3{\sf V}}(4)$ graph  with $a_{3{\sf V}}(4) = 63$ regions,
 by replacing hatpins with long narrow 3-armed ${\sf V}$'s (or Wu's).
 The new graph has $V_{\sC} = 54$ crossings and $\alpha(3{\sf V}) = 12$ arms,
 and a $\cpa$ of $9$, in agreement with \eqref{Eqcpa1}.
(Some of the arrowheads on the right have been omitted for clarity.)}
 \label{FigHatpinWu}
\end{figure}

\begin{table}[!htb]
$$
\begin{array}{c|rrrrrrrrrr}
k \backslash n & 0 & 1 & 2 & 3 & 4 & 5 & 6 & 7 & 8 & ... \\
\hline
0 &  1 &  1 &  1 &  1 &  1 &  1 &  1 &  1 &  1 &   ...  \\
1 &  1 &  2 &  4 &  7 &  11 &  16 &  22 &  29 &  37 &   ...  \\
2 &  1 &  2 &  7 &  16 &  29 &  46 &  67 &  92 &  121 &   ...  \\
3 &  1 &  3 &  14 &  34 &  63 &  101 &  148 &  204 &  269 &   ... \\
4 &  1 &  5 &  25 &  61 &  113 &  181 &  265 &  365 &  481 &   ...  \\
5 &  1 &  8 &  40 &  97 &  179 &  286 &  418 &  575 &  757 &   ...  \\
6 &  1 &  12 &  59 &  142 &  261 &  416 &  607 &  834 &  1097 &  ... \\
7 &  1 &  17 &  82 &  196 &  359 &  571 &  832 &  1142 &  1501 &  ... \\
\end{array}
$$
\caption{Table of $a_{k{\sf C}}(n)$, the maximum number of regions formed in the plane by drawing $n$ $k$-chains.  The entries are given by \eqref{EqkC2}. 
Rows 1-5 are 
  \seqnum{A000124}, \seqnum{A130883}, \seqnum{A140064}, 
  \seqnum{A080856}, \seqnum{A383465},
  columns  1 and 2 are \seqnum{A152948} and \seqnum{A386479},
  the main diagonal is \seqnum{A387525},
   and the array itself is \seqnum{A386478}.}\label{Table2}.
\end{table}

We may verify \eqref{EqkV1} directly by considering what happens to the different
 parts of the hatpin graph when the hatpins are replaced by $k{\sf V}$s.
 There are five parts in the hatpin graph, which we denote by  I (base nodes), 
 II (finite edges), III (crossing nodes),  
 IV (cells), and V (infinite edges), as shown in Fig.\ \ref{FigHatpinWu}.
 Each crossing node, for example, 
 becomes a $k \times k$  grid with $\kappa(k{\sf V}) = k^2$ vertices, $2k(k-1)$ edges, 
 and $(k-1)^2$ regions, and there are $\binom{n}{2}$ such crossing points,
 as we know from the last paragraph of~\S\ref{SecH}. By combining these counts for the five  parts
 we can verify \eqref{EqkV1}.

\section{The \texorpdfstring{$k$-chain}{k-chain}}\label{SeckC}
We define  a $k$-{\em chain} for $k \ge 1$ 
to be any continuous plane curve made up of 
$k$ piecewise linear segments. 
A line is a $1$-chain, and  a {\sf V} is a $2$-chain. 
Since it can only help in increasing the number of regions, we will assume that the outer line segments  (the first and the $k$th)  are infinite. The other segments are required to be finite. We make no restriction on the lengths of the finite segments, nor on the angles between them, and we allow self-intersections. Thus a $k$-chain for $k \ge 2$  is essentially an ordered set of $k-1$ points in the plane, with successive points joined by line segments and with semi-infinite rays attached to the first and last point. We denote the \mabs\ of all such knives
by $k{\sf C}$.  Its group is $\affineGroup$, but (as in the previous section) if $k>2$ it does not act transitively: this is also a heterotypic \mabs. (The reason is that  a $k$-chain has $2k$ degrees of freedom: $2(k-1)$ to choose the base nodes in the plane and another $2$ to describe the infinite rays.) 

A single $k$-chain may intersect itself if $k \ge 3$, producing a maximum of
$\binom{k-1}{2}$ internal regions and 
$\sigma(k{\sf C}) = \binom{k-1}{2}$ self-intersections. 
\begin{figure}[!htb]
        \begin{center}
                \begin{tikzpicture}[rotate=0]  [scale=1]

\begin{scope}[rotate=0] [scale=0.8]
\draw[ Stealth-Stealth,   thick, black] (1,5.7) -- (6,6) -- (2,6) -- (7,5.7);
\draw[ Stealth-Stealth,   thick, black] (1,5.3) -- (6.1,4.7) -- (2.1,4.9) -- (6.1,5.1) -- (1,4.55);; 
\draw[ Stealth-Stealth,   thick, black] (1,3.1) -- (6.1,4.2) -- (2,3.6) -- (6.1,3.6) -- (2.1,4.2) -- (7,3.1);
\end{scope}
                \end{tikzpicture}
        \end{center}
 \caption{Long narrow drawings of a 3-chain, a 4-chain, and a 5-chain, 
used in the proof  of~\eqref{EqkC2}. The first two are essentially unique, but there are two ways to draw the 5-chain--see Fig.\ \ref{FigMy1}.}
 \label{Fig345chains}
\end{figure}

\begin{figure}[!htb]
        \begin{center}
                \begin{tikzpicture}[rotate=0]  [scale=1]

\begin{scope}[rotate=0] [scale=1]
\draw[ Stealth-Stealth,  thick, black] (0.8,0) -- (3.2,5.30) -- (0.4,0.6) -- (4.3,6); 
\draw[ Stealth-Stealth,  thick, black] (6,0) -- (1.6,5.6) -- (5.4,1.6) -- (0,6.25);
\draw[ Stealth-Stealth,  thick, black] (0,2.9) -- (5.5,2.4) --(0.2,2.4) -- (6,2.9);
\end{scope}
                \end{tikzpicture}
        \end{center}
 \caption{Transforming a pancake graph $G_{\sf K}(3)$ to an optimal  $G_{3{\sf C}}(3)$ graph
 with $a_{3{\sf C}}(3) = 34$ regions, by replacing the knife cuts by long narrow 3-chains.}
 \label{Fig3chains3}
\end{figure}

A graph $\Gamma_{k{\sf C}}(n)$ formed by drawing $n$ $k$-chains has $V_{\sB} = n(k-1)$ base nodes and $nk$ arms, with $d_B = 2$, and $\Ei = 2n$. 
The basic equation \eqref{EqEV2}  then gives 
\beql{EqC1a}
 R = V_{\sC} + n + 1~,
 \eeq
 just as for the pancake graph (see~\eqref{EqL2}).  So again we must maximize 
$V_{\sC}$, which we can do by first maximizing 
the self-intersections in each $k$-chain, giving $n \binom{k-1}{2}$ intersections,
and then using a construction similar to that of the previous section. That is, we draw 
the individual $k$-chains as long narrow graphs (as in Fig.\ \ref{Fig345chains}),
and then substitute  them into a pancake graph if $k$ is odd 
(and the $k$-chain has  infinite edges pointing in opposite directions), 
or into a hatpin graph if $k$ is even (and 
both infinite edges of the $k$-chain point in the same direction).
See Fig.\ \ref{Fig3chains3}, which illustrates the construction for three 3-chains.
As a result we obtain $V_{\sC} = n \binom{k-1}{2} +  \binom{n}{2} k^2$, and a maximum of 
\beql{EqkC2}
R = a_{k{\sf C}}(n) =  \frac{k^2 n^2}{2} - \frac{3kn}{2} + 2n + 1~,
\eeq
regions. 
The $\cpa$ for an optimal graph is
$2 V_{\sC}/nk$, which is an integer if and only if $k \le 2$.
Table~\ref{Table2} shows the values of $a_{k{\sf C}(n)}$ for
$0 \le k \le 7$ and $0 \le n \le 8$.

The $k=2$ rows of Tables~\ref{Table1} and \ref{Table2} are of course identical,
since a 2-chain is the same as a  2-armed {\sf V}. 
The $k=1$ rows agree except that they are offset from each other by one position.  
The reason for this is 
that a 1-chain is a line, whereas a 1-armed {\sf V} is a hatpin ({\sf V}'s have a distinguished point,
lines do not.)
But the real surprise is that the $k=3$ rows of the two tables agree.
This will be explained in the next section.

\begin{figure}[!htb]
 \centerline{\includegraphics[angle=0, width=2.50in]{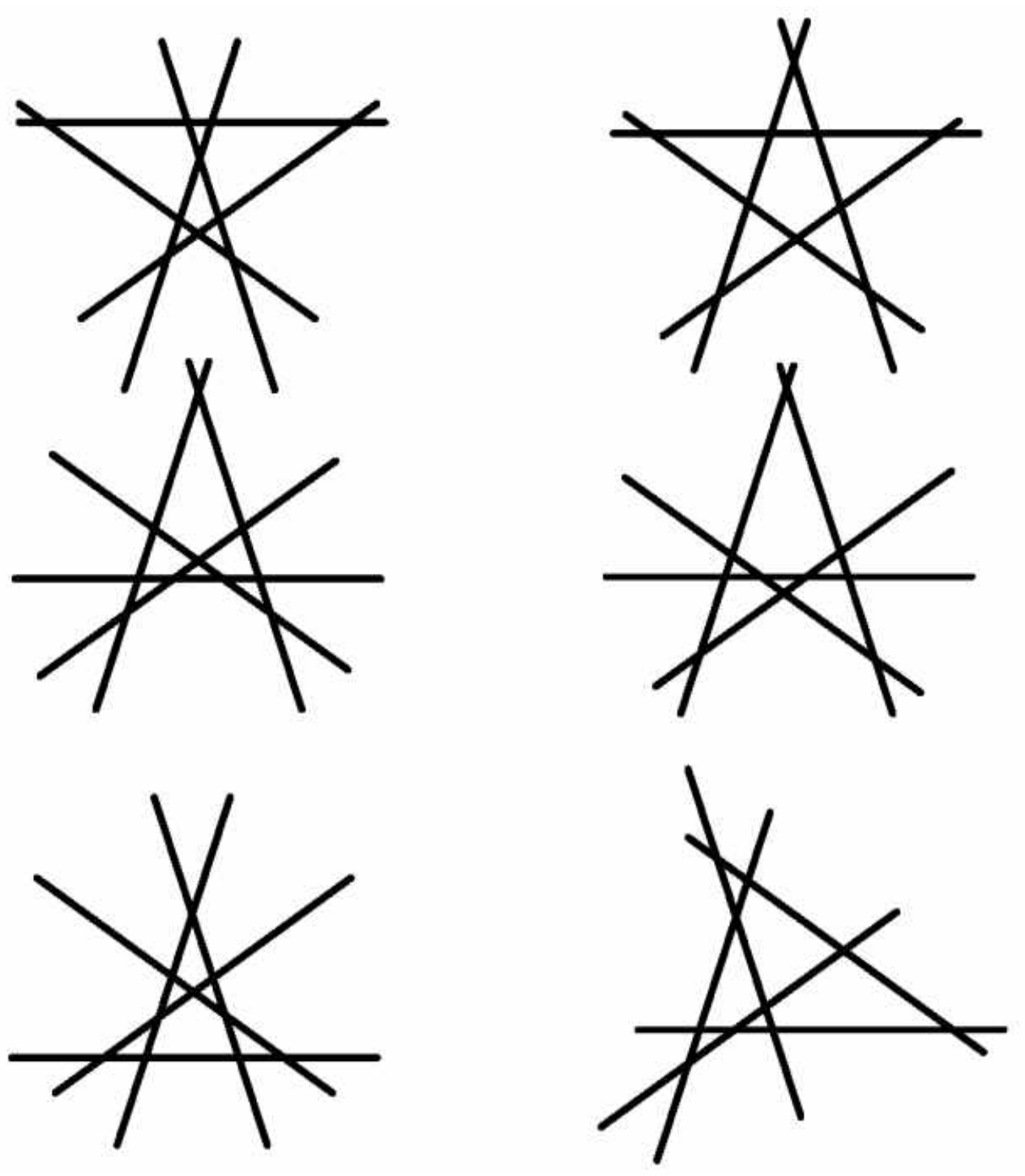}}
\caption{The six ways to draw five lines in general position in the plane [Jonathan Wild and Lawrence Reeves, from  \seqnum{A090338}]. (The arrowheads are not shown.)
Only the first two correspond to 5-chains (see Fig.\ \ref{FigMy1}).}
 \label{FigWild1}
 \end{figure}

Before leaving this section we mention another, somewhat surprising, property of Table~\ref{Table2}:
it  is very nearly symmetrical.  For example, $a_{4{\sf C}}(3)= 61$, while $a_{3{\sf C}}(4) = 63$.
This is easily explained by the formula: from \eqref{EqkC2},
$|a_{k{\sf C}}(n) - a_{n{\sf C}}(k)| = 2|k-n|$, which is exactly the difference between the numbers
of infinite regions in the two optimal graphs. Indeed, if we define $b_{k{\sf C}}(n)$ to be the maximum number of \emph{bounded} regions in any arrangement of $n$ $k$-chains, \eqref{EqkC2} implies that 
$$
b_{k{\sf C}}(n) = \frac{1}{2}(kn)^2 - \frac{3}{2}\,kn + 1,
$$
since the $2n$ edges extending to infinity create $2n$ unbounded regions. This expression is manifestly symmetric in $k$ and $n$.  So we understand the near-symmetry algebraically.

\begin{figure}[!htb]
\centering
\begin{tikzpicture}[scale=1, transform shape]   

\def\sca{2};
\def\tft{3};
\draw[ Stealth-Stealth, thick] (\sca*0,\sca*1.5) -- (\sca*1,\sca*0) -- (\sca*0,\sca*1) -- (\sca*1,\sca*1) -- (\sca*0,\sca*0) -- (\sca*1,\sca*1.5);
\draw[ Stealth-Stealth, thick] (\tft+\sca*0.4,\sca*1.5) -- (\tft+\sca*1,\sca*0) -- (\tft+\sca*0,\sca*1) -- (\tft+\sca*1,\sca*1) -- (\tft+\sca*0,\sca*0) -- (\tft+\sca*0.6,\sca*1.5);

\end{tikzpicture}
\caption{The two optimal 5-chains.}
\label{FigMy1}
\end{figure}

We can also give a geometric explanation. One possible explanation might have been that 
there is a 1-to-1-correspondence between optimal $\Gamma_{k{\sf C}}(n)$ graphs and 
optimal $\Gamma_{n{\sf C}}(k)$ graphs. But this is false even in the simplest case.
An optimal $\Gamma_{1{\sf C}}(n)$ graph consists of $n$ lines in general position in the plane, 
and the number of such graphs is known for $0 \le n \le 9$ (\seqnum{A090338}):
$$
1, 1, 1, 1, 1, 6, 43, 922, 38609, 3111341, \ldots
$$
Figure~\ref{FigWild1} shows the six ways to draw five lines in general position.
On the other hand, $\Gamma_{k{\sf C}}(1)$ is a single optimal $k$-chain 
(meaning it has  $\sigma(n) = \binom{n-1}{2}$  self-intersections) and the number of such graphs
for $0 \le k \le 5$ is $1, 1, 1, 1, 1, 2, \ldots$, the two  optimal 
five-chains being shown in Fig.\ \ref{FigMy1}.

\begin{figure}[!htb]
	\centering
	\begin{tikzpicture}[scale=1.0]
		\pgfmathsetmacro{\theta}{360/5}
		\pgfmathsetmacro{\dist}{2.5}
		\pgfmathsetmacro{\r}{sqrt((7 + 3*sqrt(5))/2)}
		\coordinate (A) at ({1}, {0});
		\coordinate (B) at ({cos(\theta)}, {sin(\theta)});
		\coordinate (C) at ({cos(2*\theta)}, {sin(2*\theta)});
		\coordinate (D) at ({cos(3*\theta)}, {sin(3*\theta)});
		\coordinate (E) at ({cos(4*\theta)}, {sin(4*\theta)});
		\draw[thick] (A) -- (B) -- (C) -- (D) -- (E) -- (A);
		
		\draw[Stealth-Stealth, thick, dashed]
		($(B)!-\dist!(A)$) -- ($(B)!1 + \dist!(A)$);
		\draw[Stealth-Stealth, thick, dashed]
		($(C)!-\dist!(B)$) -- ($(C)!1 + \dist!(B)$);
		\draw[Stealth-Stealth, thick, dashed]
		($(D)!-\dist!(C)$) -- ($(D)!1 + \dist!(C)$);
		\draw[Stealth-Stealth, thick, dashed]
		($(E)!-\dist!(D)$) -- ($(E)!1 + \dist!(D)$);
		\draw[Stealth-Stealth, thick, dashed]
		($(A)!-\dist!(E)$) -- ($(A)!1 + \dist!(E)$);
	\end{tikzpicture}
	\caption{A pentagon and its stellation.  The edges of the pentagon are replaced by lines, creating five new bounded regions and ten infinite regions.}
	\label{FigureStellation00}
\end{figure}
However, the symmetry can be explained by invoking the classical process of {\em stellation} (\cite{CDFP38}, \cite[pp.~168--173]{Cro97}, \cite[Part II]{Wen71}).
For polyhedra, stellation refers to extending the faces to create new regions, and can be carried out in stages, so that a
polyhedron can have multiple stellations. We do not need this level of detail, 
so for our planar figures, we interpret stellation
to mean that all rays and line segments are extended to infinite lines 
in a single step. 
The terminology is illustrated in  Figure~\ref{FigureStellation00}.

\begin{figure}[!htb]
	\centering
	\begin{tikzpicture}[scale = 1.3]
		\pgfmathsetmacro{\dist}{1.5}
		\coordinate (A) at (1,3.1);
		\coordinate (B) at (6.1,4.2);
		\coordinate (C) at (2,3.6);
		\coordinate (D) at (6.1,3.6);
		\coordinate (E) at (2.1,4.2);
		\coordinate (F) at (7,3.1);
		\draw[ Stealth-Stealth,   thick, black] (A) -- (B) -- (C) -- (D) -- (E) -- (F);
		\draw[-Stealth, thick, dashed]
		($(A)!0!(B)$) -- ($(A)!\dist!(B)$);
		\draw[Stealth-Stealth, thick, dashed]
		($(B)!1-\dist!(C)$) -- ($(B)!\dist!(C)$);
		\draw[Stealth-Stealth, thick, dashed]
		($(C)!1-\dist!(D)$) -- ($(C)!\dist!(D)$);
		\draw[Stealth-Stealth, thick, dashed]
		($(D)!1-\dist!(E)$) -- ($(D)!\dist!(E)$);
		\draw[Stealth-, thick, dashed]
		($(E)!1-\dist!(F)$) -- ($(E)!0!(F)$);
	\end{tikzpicture}
	\caption{When an optimal chain is stellated (here, a 5-chain), only unbounded regions are added.}
	\label{FigureStellationOptimal}
\end{figure}

So let us take an optimal configuration of $n$ $k$-chains, and stellate it. No new bounded regions can be created. 
This is because in such a configuration, any two rays or segments must already intersect: either they are consecutive elements of the same chain, in which case they share an endpoint, or else they intersect as a consequence of optimality (recall that we have already established that such configurations are possible). In either case, the stellation creates no new intersections, but only adds rays going off to infinity. This is readily seen to be the case in the optimal configurations already depicted. For example, stellating the lowermost chain in Figure~\ref{Fig345chains} produces Figure~\ref{FigureStellationOptimal}, in which only unbounded regions are created.  Figure~\ref{FigureStellation01}, on the other hand, shows a non-optimal situation.

\begin{figure}[!htb]
	\centering
	\begin{tikzpicture}
		\coordinate (A) at (-2, 2);
		\coordinate (B) at (2, -2);
		\coordinate (C) at (1.5, 0);
		\coordinate (D) at (-1, -1);
		\coordinate (E) at (2, 2);
		\draw[Stealth-Stealth, thick] (A) -- (B) -- (C) -- (D) -- (E);
		\draw[-Stealth, thick, dashed]
		($(B)!0!(C)$) -- ($(B)!2!(C)$);
	\end{tikzpicture}
	\caption{If the stellation of a chain produces new bounded regions, the chain was not optimal to start with.}
	\label{FigureStellation01}
\end{figure}
Thus, stellating an optimal configuration of $n$ $k$-chains produces a configuration of $kn$ lines, all of which intersect pairwise; and the number of bounded regions in the resulting configuration is the same as the number of bounded regions in the original chain configuration. This explains the observed symmetry.

It would  be interesting (see \S\ref{SecOP}) to know further terms in the enumeration of optimal $k$-chains, and, more generally, the numbers of graphs that achieve any of the values shown in
Tables~\ref{Table1} and \ref{Table2}.

\begin{figure}[!htb]
\centering
\begin{tikzpicture}[scale=.75, transform shape]  

\draw[ -Stealth, thick] (0.05,5.57) -- (6,1.84);
\draw[ -Stealth, thick] (0.05,5.57) -- (6,7.85);
\draw[thick,red] (1.95,4.35) -- (3.03,6.73);
\draw[ -Stealth, thick] (1.63,8.3) -- (3.3,1);
\draw[ -Stealth, thick] (1.63,8.3) -- (6,1.3);
\draw[thick,red] (2.53,6.86) -- (3.13,1.85);

\end{tikzpicture}
\caption{14 regions  can be obtained with two long-legged $A$'s.  The crossbars are shown in red.}
\label{FigA2}
\end{figure}

\section{The long-legged {\sf A}}\label{SecAG}

Our next knife is a long-legged ${\sf A}$, which  is a long-legged 2-armed {\sf V} (\S\ref{SeckV}) with a crossbar joining the two arms:
\begin{center}
\begin{tikzpicture} [scale = 0.65]
\draw[-Stealth, thick] (-4,0) -- (4,.5);
\draw[thick, red] (0,.25) -- (0.8,-.3);
\draw[-Stealth, thick] (-4,0) -- (4, -.5);
\end{tikzpicture}
\end{center}
The crossbar need not meet the arms at the same distance from the tip,
and  the angle between the arms of the {\sf V} is arbitrary. 
Equivalently, a long-legged ${\sf A}$ is a triangle where the two sides incident with a vertex have been extended to rays pointing away from that vertex. The \mabs\ ${\sf A}$ of all such knives
has symmetry group $\affineGroup$, which acts transitively.
 
Our drawings suggested that the initial values of $a_{\sf A}(n)$, the maximum 
numbers of regions achievable with $n$ ${\sf A}$'s, for $n = 1, 2$, and  3,
were 3, 14  (Fig.\ \ref{FigA2}), and 34,  matching the counts for both 3-armed Vs
and  3-chains.  This was not a coincidence!

\begin{theorem}\label{Th2} 
The maximum number of regions that can be formed in the plane
by drawing $n$ long-legged 3-armed {\sf V}'s, or $n$ 3-chains, or $n$ long-legged ${\sf A}$'s, is
$(9n^2-5n+2)/2$ (\seqnum{A140064}).
\end{theorem}

The theorem also applies to long-legged versions of the letters $E$, 
\wu, and
\begin{tikzpicture} [scale = 0.13]
\draw[thick] (0,0) -- (0,2);
\draw[thick] (-1,0) -- (-1, 2);
\draw[thick] (1,0) -- (1, 2);
\draw[thick] (-1,0) -- (1, 0); 
\end{tikzpicture}. 

\begin{figure}[!htb]
\centering
\begin{tikzpicture}[scale=1.0, transform shape]   

\draw[ -Stealth, thick] (1.5, 3.85) -- (.6, 0);
\draw[ -Stealth, thick] (1.5, 3.85) -- (2.5,0);
\draw[ -Stealth, thick] (.05,1.84) -- (3,3.5);
\draw[ -Stealth, thick] (.05,1.84) -- (3,1);
\draw[ -Stealth, thick] (3,2.3) -- (0,.8);
\draw[ -Stealth, thick] (3,2.3) -- (.2,0);
\draw[thick, red] (0.88,2.35) -- (1.2,1.54);
\draw[thick, red] (0.64,1.15) -- (1.02,.7);
\draw[ultra thick, red] (1.15,2.25) -- (2.15,1.32);
\draw[->, ultra thick, blue] (1.15,2.25) to [out=40, in=315] (1.35,3.1);
\draw[->, ultra thick, blue] (2.2,1.4) to [out=225, in=180] (2.35,0.6);

\begin{scope}[shift={(4,0)}]
\draw[ -Stealth, thick] (1.5, 3.85) -- (.6, 0);
\draw[ -Stealth, thick] (1.5, 3.85) -- (2.5,0);
\draw[ -Stealth, thick] (.05,1.84) -- (3,3.5);
\draw[ -Stealth, thick] (.05,1.84) -- (3,1);
\draw[ -Stealth, thick] (3,2.3) -- (0,.8);
\draw[ -Stealth, thick] (3,2.3) -- (.2,0);
\draw[thick, red] (0.88,2.35) -- (1.2,1.54);
\draw[thick, red] (0.64,1.15) -- (1.02,.7);
\draw[ultra thick, red] (1.35,3.1) -- (2.35,0.6);
\draw[->, ultra thick, blue] (1.35,3.1) to [out=45, in=270] (1.5,3.75);
\draw[->, ultra thick, blue] (2.35,0.6) to [out=180, in=90] (1.6,0);
\end{scope}

\begin{scope}[shift={(8,0)}]
\draw[ -Stealth, thick] (1.5, 3.85) -- (.6, 0);
\draw[ -Stealth, thick] (1.5, 3.85) -- (2.5,0);
\draw[ -Stealth, thick] (.05,1.84) -- (3,3.5);
\draw[ -Stealth, thick] (.05,1.84) -- (3,1);
\draw[ -Stealth, thick] (3,2.3) -- (0,.8);
\draw[ -Stealth, thick] (3,2.3) -- (.2,0);
\draw[thick, red] (0.88,2.35) -- (1.2,1.54);
\draw[thick, red] (0.64,1.15) -- (1.02,.7);
\draw[-Stealth, ultra thick, red] (1.5,3.85) -- (1.6,0);
\end{scope}

\end{tikzpicture}
 \caption{Converting an optimal graph with $n$ long-legged $A$'s (left)  into a graph with $n$ 3-armed $V$'s (right). The steps do not decrease the number of regions.  Reversing the last step (going from right to center) performs the inverse  operation.}
 \label{FigA2Wu}
\end{figure}

\noindent
{\em Sketch of proof.}
We have already established the theorem
for 3-armed {\sf V}'s (in~\S\ref{SeckV}) and 3-chains (in~\S\ref{SeckC}).
It remains to consider long-legged ${\sf A}$'s.

Let $G$ be an optimal arrangement of $n$ long-legged ${\sf A}$'s. Take the crossbar of any ${\sf A}$, and suppose the left end of the crossbar is closer to its tip, and the right end is further from its tip. Slide the left end of this crossbar so it is closer to the tip than any crossing node on that arm, and slide the right end of the crossbar so that it is further from the tip than any other crossing node on that arm. This is illustrated in Figure~\ref{FigA2Wu}.
We claim that after this operation has been performed on all ${\sf A}$'s, the resulting configuration is still optimal. Indeed, the legs have not been moved at all, hence all their intersections are intact. The crossbar of any ${\sf A}$ must now intersect all legs of all the other ${\sf A}$'s in the configuration, and all pairs of crossbars must intersect, as illustrated in Figure~\ref{FigureTopology}.
\begin{figure}[!htb]
	\centering
	\begin{tikzpicture}
		\draw[thick, red] (0 - 0.3, 3 - 0.9) -- (0.8, 3 - 2.4);
		\draw[thick, red] ({-2 + 3.5, 1.5 + 3.5/8}) -- ({-2 + 1, 1.5 - 1/8});
		\draw[Stealth-Stealth, thick] (-1, 0) -- (0, 3) -- (1, 0);
		\draw[Stealth - Stealth, thick] (2, 1.5-0.5) -- (-2, 1.5) -- (2, 1.5+0.5);
	\end{tikzpicture}
	\caption{After sliding, the crossbars join alternate sides of the quadrilateral formed by the legs of the {\sf A}s and so must intersect.}
	\label{FigureTopology}
\end{figure}

Finally, we can move the left end of the crossbar so it actually coincides with the tip and the right end so it is cut free from its arm. This ${\sf A}$ has now become a 3-armed {\sf V}, and
since $G$ was optimal, the number of regions has
remained  unchanged throughout the process.
By repeating this operation for all the ${\sf A}$'s in $G$  we turn the ${\sf A}$ graph  
into a 3-armed {\sf V}  graph. The converse operation is equally simple. 
\begin{figure}[!htb]
\centering
\begin{tikzpicture}[scale=0.22, transform shape]   

\draw[Stealth-, thick, black] (2,1) -- (4,8);
\draw[Stealth-,thick, black] (10,1) -- (8.6,5);
\draw[thick, black] (4,8) -- (5.8,14.4);
\draw[thick, red] (4,8) -- (8.6,5);
\draw[thick, black] (5.8,14.4) -- (8.6,5);

\begin{scope}[shift={(10,0)}]
\draw[Stealth-Stealth, thick, black] (2,1) -- (9.5,27);
\draw[Stealth-Stealth, thick, black] (10,1) -- (1.5,27);
\draw[thick, red] (4,8) -- (8.6,5);
\draw[thick, black] (5.8,14.4) -- (8.6,5);
\draw[->, ultra thick, blue] (3.3,5.5) to [out=135, in=190] (8,22);
\draw[->, ultra thick, blue] (9.3,2.6) to [out=20, in=20] (4.5,18);
\end{scope}

\begin{scope}[shift={(20,0)}]
\draw[-Stealth, thick, black] (4,8) -- (9.5,27);
\draw[-Stealth, thick, black] (8.6,5) -- (1.5,27);
\draw[thick, red] (4,8) -- (8.6,5);
\end{scope}

\end{tikzpicture}
\caption{Converting a long-legged $A$ (left) into a 3-chain (right), }
\label{FigA23chain}
\end{figure}
On the other hand, to convert the long-legged ${\sf A}$ graph $G$ into an arrangement of $3$-chains, we simply apply the transformation shown in Fig.\ \ref{FigA23chain} to each ${\sf A}$. This preserves the number of regions (we invite the reader to apply it to Fig.\ \ref{FigA2} and check that the result still has 14 regions).

\section{The wynn}\label{SecWynn}
In this section we consider the Runic (and Middle English) letter wynn (\runewynn). 
One reason this \mabs\ of knives is interesting is that it provides a counterexample 
to a conjecture we once made, that for any \mabs\ in the affine genus it would be
possible to achieve equality in Equation~\eqref{EqVx}. 
This would have been a powerful result, for it
would have meant that the maximal number of crossings between any number of copies of an affine shape
was limited solely by the number of self-crossings and the number of crossings between two copies,
with no additional restrictions required for three or more copies. The wynn \mabs\
will show  that our belief was too optimistic. 
\begin{figure}[!htb]
	\centering
	\begin{tikzpicture} [scale=1]
		\draw[thick, -Stealth] (0, 0) -- (0, -2);
		\draw[thick] (0, 0) -- (1, -.5) -- (0, -1);
	\end{tikzpicture}
	\caption{A wynn.}
	\label{FigureWynnDefinition}
\end{figure}

We define a wynn to be a (non-degenerate) triangle, one of whose sides has been 
extended to an infinite ray (Fig.\ \ref{FigureWynnDefinition}).  One instinctively draws an isosceles triangle, 
but the triangle may be equilateral, isosceles, or scalene.
Any two such figures are equivalent under the affine group $\affineGroup$, since
as already mentioned in \S\ref{SecLLL},  this group is 3-transitive. 
\begin{figure}[!htb]
	\centering
	\begin{tikzpicture}[scale=2]
		\coordinate (O) at (0, 0.25);
		\begin{scope}[rotate around={7:(O)}]
			\coordinate (A) at (3, 0);
			\coordinate (B) at (-1, 0);
			\coordinate (C) at (0, 1);
			\coordinate (E) at (2, -1);
			
			\path[name path=AB] (A) -- (B);
			\path[name path=CE] (C) -- (E);
			
			\path[name intersections={of=AB and CE, by=D}];
			
			\draw[thick, Stealth-] (A) -- (B) -- (C) -- (D);
		\end{scope}
		\begin{scope}[rotate around={-7:(O)}]
			\coordinate (A1) at (3, 1/2);
			\coordinate (B1) at (-1, 1/2);
			\coordinate (C1) at (0, -1/2);
			\coordinate (E1) at (2, 3/2);
			
			\path[name path=AB1] (A1) -- (B1);
			\path[name path=CE1] (C1) -- (E1);
			
			\path[name intersections={of=AB1 and CE1, by=D1}];
			
			\draw[thick, Stealth-] (A1) -- (B1) -- (C1) -- (D1);
		\end{scope}
		
	\end{tikzpicture}
	\caption{Two wynns can intersect in 7 points and divide the plane into 10 regions.}
	\label{FigureWynnWynn}
\end{figure}
Two wynns can intersect in at most seven points, and this can only happen if they are related by an improper affine map, that is, one that reverses orientations (Fig.\ \ref{FigureWynnWynn}).
This immediately implies
that three copies of the wynn cannot be placed so that each pair has seven intersections: there are only two orientations! The best we can do, given $n$ wynns, is to split them as evenly as possible into two groups, and arrange them so that two wynns in different groups intersect in seven points, while two wynns in the same group intersect in six points. This is possible, because the six intersections between two wynns with the same orientation can be achieved by a small rotation of the triangle, as in Figure~\ref{FigureWynnRotation}.

\begin{figure}[!htb]
	\centering
	\begin{tikzpicture}[scale=2]
		\coordinate (O) at (0, 0.25);
		\begin{scope}[rotate around={0:(O)}]
			\coordinate (A) at (3, 0);
			\coordinate (B) at (-1, 0);
			\coordinate (C) at (0, 1);
			\coordinate (E) at (2, -1);
			
			\path[name path=AB] (A) -- (B);
			\path[name path=CE] (C) -- (E);
			
			\path[name intersections={of=AB and CE, by=D}];
			
			\draw[thick, Stealth-] (A) -- (B) -- (C) -- (D);
		\end{scope}
		\begin{scope}[rotate around={5:(O)}]
			\coordinate (A1) at (3, 0);
			\coordinate (B1) at (-1, 0);
			\coordinate (C1) at (0, 1);
			\coordinate (E1) at (2, -1);
			
			\path[name path=AB1] (A1) -- (B1);
			\path[name path=CE1] (C1) -- (E1);
			
			\path[name intersections={of=AB1 and CE1, by=D1}];
			\draw[thick, Stealth-] (A1) -- (B1) -- (C1) -- (D1);
		\end{scope}
	\end{tikzpicture}
	
	\caption{Wynns with the same orientation can intersect in six points.}
	\label{FigureWynnRotation}
\end{figure}

The resulting configuration is optimal, since  we have the largest number of pairs with seven intersection points, and all remaining pair have six intersection points. Thus, with $m = \lfloor n/2 \rfloor$ and $m' = n - m$, we have
$$
V_{\sC} = 7mm' + 6\left(\binom{m}{2} + \binom{m'}{2}\right).
$$
The number of regions is then given by $R = V_{\sC} + n + 1$, and so
\beql{EqWynn}
a_{\runewynn}(n) = 3n^2 - 2n + 1 + \left\lfloor \frac{n^2}{4} \right\rfloor  \quad  (\seqnum{A393448}) ~.
\eeq
This argument---finding optimal configurations for a small number of shapes and then doubling up the shapes---will be used again for  the constrained long-legged $\overline{\sf T}$.

\section{The constrained long-legged
	\texorpdfstring{$\overline{{\sf Z}}$}{\sf Z},
	\texorpdfstring{$\overline{{\sf W}}$}{\sf W} and
	\texorpdfstring{$\overline{{\sf M}}$}{\sf M}}\label{SecZMW}
Starting in this section we consider several \mabp\ of constrained long-legged knives for which 
the symmetry group is  the smaller group $\simGroup$ of \emph{similarities} of the plane: rigid motions and scalings. Now the notions of angle and parallelism become meaningful.
If we were still working with the affine group, the shapes we are about to consider would simply be special cases of chains. The formula \eqref{EqEV2} for the number of regions is still valid, as it depends only on the graph formed by the figures. The role of the group is just to help in defining the class of figures to be considered.

For the long-legged $\overline{{\sf Z}}$ knife mentioned  in the Introduction
(the zig-zag curve of \cite{GKP94}), the basic equation \eqref{EqEV2}  gives
\eqref{EqL2} again. 
To maximize 
$V_{\sC}$ we substitute  $n$ long narrow copies of $\overline{{\sf Z}}$ into the pancake graph,
which results in a maximum of $R = a_{\overline{{\sf Z}}}(n) = (9n^2-7n+2)/2$ regions
(see \seqnum{A117625}). This is exactly the construction used in \cite{GKP94}
(although they do not use our counting argument to
prove it is optimal).

A natural sequel is worth mentioning. After drawing the first three arms of the $\overline{{\sf Z}}$,
we may truncate the third arm and add a fourth (infinite)  arm
that is  parallel to the second arm, obtaining this long-legged $\overline{{\sf W}}$ (or $\overline{{\sf M}}$) knife:
\begin{center}
	\begin{tikzpicture} [scale = 0.50]
		\draw[ Stealth-, thick] (-4,.5) -- (.5,.5);
		\draw[thick] (.5,.5) -- (-.5,-.5);
		\draw[thick] (-.5,-.5) -- (1, -.5);
		\draw[-Stealth, thick] (1, -.5) -- (-1, -2.5);
	\end{tikzpicture}
\end{center}
Our standard analysis 
(this time substituting into a hatpin graph) shows that the maximum
number of regions is $a_{\overline{{\sf W}}}(n) = 8n^2-7n+1$ (\seqnum{A125201}).
Figure \ref{FigLLM2} illustrates $a_{\overline{{\sf W}}}(2) = 19$. The parallelism constraints still allow us to make this letter arbitrarily long and thin, so that the maximum is indeed achievable. 

\begin{figure}[!htb]
\centering
\begin{tikzpicture}[scale=.5, transform shape]   

\draw[ Stealth-Stealth, thick, blue] (0.8,3.9) -- (6.7, 3.9) -- (1.6,2.8) -- (6.7, 2.8) -- (0.8,1.5);
\draw[ Stealth-Stealth, thick, black] (3,4.9) -- (4.5,0.5) -- (4.5,4.4) -- (5.8,0.6) -- (5.8,4.9);

\end{tikzpicture}
\caption{Two long-legged $\overline{{\sf W}}$'s divide the plane into 19 regions.}
\label{FigLLM2}
\end{figure}

\section{The  constrained long-legged  
\texorpdfstring{$\overline{{\sf L}}$}{\sf L},
\texorpdfstring{$\overline{{\sf X}}$}{\sf X},
\texorpdfstring{$\overline{{\sf H}}$}{\sf H},
and 
\texorpdfstring{$\overline{{\phi}}$}{phi}}\label{SecLX}
In this section we discuss another four shapes from the similarity genus, 
$\overline{\sf L}$,  $\overline{{\sf X}}$, $\overline{{\sf H}}$, and $\overline{{\phi}}$, all of which we are able to solve.

\subsection{The constrained long-legged  \texorpdfstring{$\overline{\sf L}$}{\sf L}} \label{SecSLLL}
With no restriction on the angle, a long-legged ${\sf L}$ would be
indistinguishable from a long-legged {\sf V}. 
We define a {\em constrained long-legged} $\overline{{\sf L}}$  to consist of a distinguished point, its center,  
 from which two perpendicular rays emanate.  We denote the \mabs\ of all such knives
by ${\sf L}$.  The group  $\simGroup$ acts  transitively and as holotype we may take
the center $(0,0)$ and the rays $\{y=0, x \ge 0\}$  and $\{x=0, y \ge 0\}$. 
 Two $\overline{{\sf L}}$'s can intersect in at most $\kappa({\overline{\sf L}}) = 3$ points, 
so $V_{\sC} \le 3\binom{n}{2}$ and then \eqref{EqEV2} gives
\beql{EqLLL1}
R \le 3\binom{n}{2} + n + 1 = (3n^2 - n + 2)/2 \quad (\seqnum{A143689}) ~.
\eeq

\begin{figure}[!htb]
	\centering
	\begin{tikzpicture}[scale=2]
		\draw (0, 0) circle[radius = 1];
		\anglerays{0}
		\anglerays{18}
		\anglerays{36}
		\anglerays{54}
		\anglerays{72}	
	\end{tikzpicture}
	\caption{The construction of $n$ constrained long-legged $\overline{{\sf L}}$'s, showing the case  $n=5$. There are $36$ regions.}
	\label{FigSL5}
\end{figure}

Equality in \eqref{EqLLL1} can be achieved by modifying the standard 
 construction for drawing  $n$ squares with the maximum
number of regions (see \seqnum{A069894}). This construction 
starts by marking $4n$ equally-spaced points $P_0, \ldots, P_{4n-1}$ around 
the full circumference of a circle. To get an optimal arrangement of $n$ squares,
we would take the squares with vertices $\{P_i, P_{i+n}, P_{i+2n}, P_{i+3n}\}$
for $i = 0, 1, \ldots, n-1$. To get an optimal arrangement of $n$ $\overline{{\sf L}}$'s,
we  take the $\overline{{\sf L}}$'s formed by the top left corners of these squares.  That is,  we draw
an $\overline{{\sf L}}$ centered at $P_i$ with arms $P_i - P_{i+n}$
and $P_i -  P_{i+3n}$, for  $i=0,\ldots,n-1$ .  
Figure~\ref{FigSL5} shows the construction when $n=5$.
In general the  graph contains $2n$ infinite regions,  there are $n$ rows of $n-1$ cells each, 
and a further triangle of $(n-1)(n-2)/2$ cells, for a total of
$R = a_{\overline{{\sf L}}}(n) = (3n^2 - n +2)/2$ regions. 
For $n \ge 0$ the values are $1, 2, 6, 13, 36, 52, 71, \ldots$
 (\seqnum{A143689}). The OEIS entry makes no mention of this geometric application, 
 so this result may be new.
 
 There are $V = n + 3 \binom{n}{2} = n(3n-1)/2 = R-1$ vertices, which is a pentagonal
 number (\seqnum{A000326}).
 The presence of  pentagonal numbers can be
 explained as follows.  The familiar picture of the pentagonal numbers
  (\cite[p.\ 380]{GKP94}, for example)
 shows  a pentagonal array of dots divided into successive shells of $1, 4, 7, 10, 13, \ldots$
 dots.  Careful examination of the grid of points 
 formed by the above construction (as in Fig.\ \ref{FigSL5})
 reveals exactly the same shells of points.

\begin{figure}[!htb]
\centering
\begin{tikzpicture}[scale=.75, transform shape]   
\newcommand*{\TickSize}{2pt}%

\draw [axis line style] (-1,0) -- (9,0);
\draw [axis line style] (0,-3) -- (0,5);

\foreach \x in {0,1,...,9} {%
    \draw ($(\x,0) + (0,-\TickSize)$) -- ($(\x,0) + (0,\TickSize)$)
        node [below] {$\x$};
}

\foreach \y in {-3,-2,...,-1} {%
    \draw ($(0,\y) + (-\TickSize,0)$) -- ($(0,\y) + (\TickSize,0)$)
        node [left] {$\y$};
}
\foreach \y in {1,2,...,5} {%
    \draw ($(0,\y) + (-\TickSize,0)$) -- ($(0,\y) + (\TickSize,0)$)
        node [left] {$\y$};
}
\draw[Stealth-Stealth, ultra thick, blue] (-1,0) -- (9,0);
\draw[Stealth-Stealth, ultra thick, blue] (-1,-.45) -- (9,1.95);
\draw[Stealth-Stealth, ultra thick, blue] (-1,-1.67) -- (9,3.90);
\draw[thick] (0,-2.80) -- (8,4.90);
\draw[thick] (2,-3.15) -- (7,5.00);
\draw[thick] (4.15,-3.00) -- (6.4,5.00);
\draw[Stealth-Stealth, ultra thick, blue] (6,-3.00) -- (6,5.00);
\draw[Stealth-Stealth, ultra thick, blue] (7.75,-3.00) -- (5.6,5.00);
\draw[Stealth-Stealth, ultra thick, blue] (9.70,-3.00) -- (5.15,5.00);
\draw[red, ultra thick, fill] (6,0) circle [radius=.075];
\draw[red, ultra thick, fill] (6.54,1.38) circle [radius=.075];
\draw[red, ultra thick, fill] (6.56,2.57) circle [radius=.075];
                \end{tikzpicture}
\caption{Dividing the plane into the maximum number of pieces
using $n$ $\overline{\sf X}$'s (or coordinate frames).
For $n=3$, 
draw a grid of  lines $\lambda_t$ defined by \eqref{EqPancake} with 
$\theta = \pi/(4n) = \pi/12$ and $t=0,\ldots,3n-1=8$.
The three $\overline{\sf X}$'s  (centered at the red dots, and with  blue axes) are
$\lambda_0 \cup \lambda_6$,
$\lambda_1 \cup \lambda_7$, and 
$\lambda_2 \cup \lambda_8$.
An enlargement of the figure will show that,
although it is not clear from this illustration,
there is no point where any three of the lines intersect.}
\label{FigcXfromK}
\end{figure}
 
\subsection{The constrained long-legged  \texorpdfstring{$\overline{\sf X}$}{\sf X}} \label{SeccX} 
A {\em constrained long-legged} $\overline{\sf X}$    consists of a distinguished point, its center, 
from which two perpendicular lines  (its arms, or  axes) emanate.
The group $\simGroup$ acts transitively on this \mabs, and as holotype we may take
the center $(0,0)$ and the lines $x=0$ and $y=0$.
  
A graph $\Gamma_{\overline{\sf X}}(n)$  with $n$ $\overline{\sf X}$'s is also a pancake graph $\Gamma_{\sf K}(2n)$ with $2n$ lines, 
 so \eqref{Eq1} gives 
 \beql{EqcX}
 a_{\overline{\sf X}}(n) \le a_{\sf K}(2n) = 2n^2+n+1 \quad (\seqnum{A084849})~.
 \eeq
Two $\overline{\sf X}$'s can intersect in at most $\kappa(\overline{\sf X}) = 4$ points,
so $V_{\sC} \le 4\binom{n}{2}$, and 
then \eqref{EqEV2} also gives  $R \le 2n^2+n+1$.

We can indeed achieve equality in \eqref{EqcX} by modifying the construction 
of pancake graphs given in \S\ref{SecPG}. To construct a $\Gamma_{\overline{\sf X}}(n)$
containing $n$ copies of $\overline{\sf X}$,  set $\theta = \pi/(4n)$, and let $\lambda_t$
denote the line \eqref{EqPancake}. Draw the lines $\lambda_t$ for
$t=0,1, \ldots, 3n-1$ (now we need more of these  lines  than 
 when we  constructed the pancake graph). The $t$-th copy of $\overline{\sf X}$ is 
 formed from the lines $\lambda_t$ and $\lambda_{t+2n}$ 
 for $t=0, \ldots,  n-1$. Figure~\ref{FigcXfromK}
illustrates the construction for $n=3$.
Of all our constructions, this one seems to be the most likely to have applications,
since the $\overline{\sf X}$ suggests the cross-hairs of an optical instrument.

\begin{figure}[!htb]
\begin{center}
\begin{tikzpicture} [scale = 1.6] 

\def\s{0}
\draw[thick] (0+\s,1) -- (0+\s,-1); 
\draw[thick] (1+\s,1) -- (1+\s,-1); 
\draw[ultra thick] (0+\s,0) -- (1+\s,0); 

\def\s{2}
\draw[thick] (0+\s,0) -- (-.3+\s,1); 
\draw[thick] (0+\s,0) -- (+.3+\s,1); 
\draw[thick] (1+\s,0) -- (+.7+\s,1); 
\draw[thick] (1+\s,0) -- (1.3+\s,1); 
\draw[ultra thick] (0+\s,0) -- (1+\s,0); 

\def\s{4}
\draw[thick] (-0.3+\s,1) -- (1+\s,0); 
\draw[thick] (0+\s,1) -- (1+\s,0); 
\draw[thick] (1.3+\s,1) -- (0+\s,0); 
\draw[thick] (1+\s,1) -- (0+\s,0); 
\draw[ultra thick] (0+\s,0) -- (1+\s,0); 

\def\s{6}
\draw[thick] (-0.3+\s,1) -- (1+\s,0); 
\draw[thick] (0+\s,1) -- (1+\s,0); 
\draw[thick] (1.3+\s,1) -- (0+\s,0); 
\draw[thick] (1+\s,1) -- (0+\s,0); 
\draw[ultra thick] (0+\s,1) -- (1+\s,1); 

\end{tikzpicture}
\end{center}
\caption{An unconstrained long-legged {\sf H} is the same as two {\sf V}'s  joined by a line segment, which can be rearranged to have four self-crossings, and then (invertibly) 
transformed into a 5-chain by ``raising the bar'' (the heavy line).}
 \label{FigH25}
\end{figure}

\subsection{The constrained long-legged 
       \texorpdfstring{$\overline{\sf H}$}{\sf H}
and \texorpdfstring{$\overline{\phi}$}{phi}} \label{SeccH} 
Our rules for long-legged letters tell us that an (unconstrained) long-legged {\sf H} is the same as 
two long legged {\sf V}'s  with their tips joined by a line segment (see Fig.~\ref{FigH25}).
Then $\sigma({\sf H}) = 4$, $\kappa({\sf H}) = 25$, and with $n$ {\sf H}'s, \eqref{EqEV} 
tells us that $a_{\sf H}(n) \le (25 n^2 - 11 n + 2)/2$. This upper bound matches the counts for a $5$-chain, \eqref{EqkC2} (\seqnum{A383465}), 
and equality is easily achieved since the invertible transformation shown in Fig.~\ref{FigH25}
maps an optimal graph with $n$ {\sf H}'s into one with $n$ $5$-chains.

\begin{figure}[!htb]
\begin{center}
\begin{tikzpicture}
	\coordinate (A) at (0,0);
	\coordinate (B) at (1.5,0);
	\coordinate (B1) at ($(B)!2!90:(A)$);
	\coordinate (B2) at ($(B)!3!-90:(A)$);
	\coordinate (A1) at ($(A)!3!90:(B)$);
	\coordinate (A2) at ($(A)!2!-90:(B)$);	
	\draw[thick, name path=AB] (A) -- (B);
	\draw[thick, name path=AL] (A1) -- (A2);
	\draw[thick, name path=AR] (B1) -- (B2);
	\coordinate (C) at (0.75, -0.25);
	\coordinate (D) at ($ (C) + (40:1.5) $);
	\coordinate (D1) at ($(D)!2!90:(C)$);
	\coordinate (D2) at ($(D)!3!-90:(C)$);
	\coordinate (C1) at ($(C)!3!90:(D)$);
	\coordinate (C2) at ($(C)!2!-90:(D)$);
	\draw[thick, red, name path=CD] (C) -- (D);
	\draw[thick, red, name path=CL] (C1) -- (C2);
	\draw[thick, red, name path=CR] (D1) -- (D2);
	\path[name intersections={of=AB and CD, by=I1}];
	\fill[red] (I1) circle (2pt);
	\path[name intersections={of=AB and CL, by=I2}];
	\fill[red] (I2) circle (2pt);
	\path[name intersections={of=AR and CD, by=I3}];
	\fill[red] (I3) circle (2pt);
	\path[name intersections={of=AR and CR, by=I4}];
	\fill[red] (I4) circle (2pt);
	\path[name intersections={of=AL and CR, by=I5}];
	\fill[red] (I5) circle (2pt);
	\path[name intersections={of=AL and CL, by=I6}];
	\fill[red] (I6) circle (2pt);
	\path[name intersections={of=AR and CL, by=I7}];
	\fill[red] (I7) circle (2pt);
\end{tikzpicture}
\end{center}
\caption{One of the three inequivalent ways in which two constrained $\overline{\sf H}$'s
can interesct in seven points.}
 \label{FigcH2}
\end{figure}

\tikzset{
	Hshape/.pic={
		\draw[thick] (-1, 4) -- (-1,-3);
		\draw[thick] (1, 4) -- (1,-3);
		\draw[thick] (-1,0) -- (1,0);
	}
}
\begin{figure}[!htb]
	\begin{center}		
		\begin{tikzpicture}[scale=2]
			\begin{scope} 
				\clip (-1.5,-1) rectangle (2, 3);
				\begin{scope}
					\pic[transform shape] at (0,0) {Hshape};
					\pic[transform shape, rotate=10] at (0.1, 0) {Hshape};
					\pic[transform shape, rotate=30] at (0.3, 0) {Hshape};
				\end{scope}
			\end{scope}
		\end{tikzpicture}
	\end{center}
	\caption{Construction of $n$ $\overline{\sf H}$'s 
		with each pair intersecting in 7 points (here $n=3$).}
	\label{FigcH3J}
\end{figure}

We define a {\em constrained long-legged}  $\overline{\sf H}$ to consist of a pair of parallel lines
 joined by a perpendicular line segment (the crossbar).  The group $\simGroup$ acts 
 transitively on the \mabs\ of all such knives.
Two $\overline{\sf H}$'s can meet in 
a maximum of $\kappa(\overline{\sf H}) = 7$
points. This can be done in three inequivalent ways, but we will only make 
use of one of them,  the arrangement shown in Fig.\  \ref{FigcH2}.

It is straightforward to arrange $n$ $\overline{\sf H}$'s so that each pair intersect in seven points
in this way, 
by first drawing the $n$ crossbars so that each one crosses  the others at a small angle, similar to the 
arrangement of lines in Fig.\ \ref{FigKnife2}. By making these crossbars have slowly increasing and incommensurable lengths, we ensure that when we add the parallel lines to the crossbars,
there are no triple intersections in the resulting graph. 
Figure \ref{FigcH3J} illustrates the construction when $n=3$.

From \eqref{EqEV}, the maximum number of regions with $n$ $\overline{\sf H}$'s is then
 \beql{EqcH}
 a_{\overline{\sf H}}(n) = \frac{7n^2-n+2}{2} \quad (\seqnum{A140063})~.
 \eeq
  with initial values $1, 4, 14, 31, 55, \ldots$.
We were surprised to find that this sequence was already in the OEIS, in a 2008 entry.
 This entry formerly contained a link to a webpage \cite{Kin09}  of David Kinsella, 
 where it was stated that this sequence gives 
 the maximum number of regions that the plane can be divided into by drawing $n$ circles and $n$ lines.
 This webpage now seems to be lost, but the Internet Archive's {\em Wayback Machine} 
 has preserved enough of this document (see \cite{Kin09}: the figures are missing) 
 to make it clear that this was an empirical observation based on examination of the 
 constructions for $n \le 4$ pairs of circles and lines.
 
 \begin{figure}[!htb]
\begin{center}
\begin{tikzpicture} [scale = 1.5] 


\def\s{0}
\draw[thick] (0+\s,1) -- (0+\s,-1); 
\draw[thick] (1+\s,1) -- (1+\s,-1); 
\draw[thick] (0+\s,0) -- (1+\s,0); 

\def\s{3.8}
\draw[thick] (0+\s,1) -- (0+\s,-1); 
\draw[thick] (1+\s,1) -- (1+\s,-1); 
\draw[thick ] (0+\s,0) -- (1+\s,0); 
\draw[thick ] (-1.2+\s,0) -- (0+\s,0); 
\draw[thick ] (1+\s,0) -- (2.2+\s,0); 
\draw[thick] (0.5+\s,0)  circle  (0.5cm);
\draw[->, ultra thick, blue] (0+\s,-.5) to [out=180, in=270] (-.5+\s,0);
\draw[->, ultra thick, blue] (0+\s,.8) to [out=0, in=90] (.5+\s,0.5);
\draw[->, ultra thick, blue] (1+\s,.5) to [out=0, in=90] (1.5+\s,0);
\draw[->, ultra thick, blue] (1+\s,-.8) to [out=180, in=300] (.5+\s,-0.5);

\def\s{7.6}
\draw[thick ] (0+\s,0) -- (1+\s,0); 
\draw[thick ] (-1.3+\s,0) -- (0+\s,0); 
\draw[thick ] (1+\s,0) -- (2.3+\s,0); 
\draw[thick] (0.5+\s,0)  circle  (0.5cm);

\end{tikzpicture}
\end{center}
\caption{An invertible map from $\overline{\sf H}$'s to $\overline{\phi}$'s.}
\label{FigH2phi}
\end{figure}

 A stronger version of this claim can be obtained from our analysis of the $\overline{\sf H}$ shape.
 Define a {\em constrained}  $\overline{\phi}$ to consist of a circle with a line through its center.
 Two such $\overline{\phi}$'s can intersect in a maximum of seven points. An optimal
 graph containing $n$ $\overline{\phi}$'s, each intersecting all of the others in seven points, can be 
 obtained from an optimal graph with $n$ $\overline{\sf H}$'s by applying the invertible transformation shown in Fig.\ \ref{FigH2phi}. This map does not preserve the number of intersections between the figures in general, but an elementary check shows that it \emph{does} map two shapes with seven intersections to two shapes with seven intersections. It follows that \eqref{EqcH} applies also to $\overline{\phi}$'s. Thus, not only is Kinsella's claim true, but it is possible to achieve the maximum number of regions under the additional constraint that each line must pass through the center of a circle.

 \subsection{Steiner's theorem:
 \texorpdfstring{$m$ circles and $n$ lines}{m circles and n lines}} \label{SecSteiner}
Unsurprisingly, the combination of  $n$ circles and $n$ lines had already been 
discussed in the literature. In fact,  Steiner \cite{Ste26} had already  shown  that the maximum
number of regions that the plane can be divided into by making $m$ independent
cuts with circular knives from the \mabs\ ${\sf O}$ and at the same time making $n$ independent
cuts with straight knives from ${\sf K}$ is given by
\begin{displaymath}
        a_{{\sf St}}(m,n) = \begin{cases}
        m^2-m+2mn+(n^2+n+2)/2,  & \text{ if $n>0$ or  $m=n=0$;}\\
       m^2+m+2, & \text{ if $n=0$ and $m>0$ . }\\           \end{cases}
 \end{displaymath}
The diagonal terms $a_{{\sf St}}(n,n)$ agree with \eqref{EqcH}. 
The OEIS entry \seqnum{A393442} contains the array
$a_{{\sf St}}(m,n)$  itself (analogous to Tables~\ref{Table1} and \ref{Table2}), the first three columns
of which  are essentially \seqnum{A014206}, and the first two rows are \seqnum{A000124} and \seqnum{A034856}.

We remark that this \mabs\ of knives does not fit into our usual classification into genera, since now we are choosing a single knife from the very large \mabs\ ${\sf St}(m,n)$ rather than making $n$ independent choices from a smaller \mabs.

It is also worth remarking that Steiner's article \cite{Ste26} appeared exactly 200 years ago
 in the first volume  of what became one of the most famous of all mathematical journals, the 
{\em Journal f\"{u}r die reine und angewandte Mathematik}, often simply called {\em Crelle's Journal}.
Steiner's proof  of his formula essentially uses a case-by-case analysis, 
considering all possible ways in which $m$ circles can be divided into families of parallel circles 
and all possible ways in which  $n$ lines can be divided into families of parallel lines. The article  occupies 15 pages of {\em Crelle's Journal}, and we will not present the proof here. Steiner's proof is also discussed by Wetzel \cite{Wet78}.

The many citations of Steiner's paper in the literature contain a great deal of information 
about dissections by lines and circles and refinements (going beyond simply maximizing the number of pieces, for example), generalizations  (e.g., to higher dimensions), and applications.
At present, Google Scholar lists over 130 citations. 
See for example \cite{Ale81, Bar19, Ede86, Mil82,  Pac95,  Rin56, Wet78}.

\section{The  constrained long-legged  \texorpdfstring{$\overline{{\sf T}}$}{\sf T}   and 
\texorpdfstring{$\overline{{\sf A}}$}{{\sf A}}}\label{SecTA}
Up to this point in our investigations we were able to analyze the shapes by hand,
using pencil and paper and traditional drawing instruments. But the next two \mabp\ in the similarity genus,
$\overline{\sf T}$ and $\overline{{\sf A}}$, appear to be more difficult and we invoked the aid of the computer. In the case of $\overline{\sf T}$, we find an exact formula, which, however, is no longer a polynomial in $n$ but a quasipolynomial. The case of $\overline{{\sf A}}$ is harder and here we are only able to determine the numbers up to $n = 3$. Beyond that point, we only have lower and upper bounds.

\begin{figure}[!htb]
\centering
\begin{tikzpicture}[scale=.15, transform shape]   

\draw[Stealth-Stealth, thick] (0,5) -- (34,5); 
\draw[-Stealth, thick] (15,5) -- (15,26); 
\draw[Stealth-Stealth, thick] (11,26) -- (26.3,0);
\draw[-Stealth, thick] (18.9,12.7) -- (0,2);
\draw[Stealth-Stealth, thick] (5.8,0) -- (21,26);
\draw[-Stealth, thick] (13.5,13.2) -- (34,1.5);

                \end{tikzpicture}
\caption{Three constrained long-legged $\overline{\sf T}$'s can divide 
the plane into 19 regions.  It is impossible to add a fourth 
 $\overline{\sf T}$
that meets each of the first three in 4 points.}
 \label{FigT19}
\end{figure}

\begin{figure}[!htb]
\centering
\begin{tikzpicture}[scale=.25, transform shape]   

\draw[Stealth-Stealth, thick] (0,5) -- (34,5); 
\draw[-Stealth, thick] (15,5) -- (15,26); 
\draw[Stealth-Stealth, thick] (11,26) -- (26.3,0);
\draw[-Stealth, thick] (18.9,12.7) -- (0,2);
\draw[Stealth-Stealth, thick] (5.8,0) -- (21,26);
\draw[-Stealth, thick] (13.5,13.2) -- (34,1.5);
\draw[Stealth-Stealth, thick, red] (7.0,0) -- (12.6,26);
\draw[-Stealth, thick, red] (10.00,13.90) -- (34,7.3);

                \end{tikzpicture}
\caption{Four constrained long-legged $\overline{\sf T}$'s can divide 
the plane into 32 regions. The first three are the same as in Fig.\ \ref{FigT19},
and the fourth (which meets one of the first three in only three points, 
not four) is drawn in red.}
 \label{FigT32}
\end{figure}

\subsection{The constrained long-legged  \texorpdfstring{$\overline{{\sf T}}$}{\sf T}}\label{SeccT}
A {\em constrained long-legged} $\overline{\sf T}$  consists of a distinguished point, its center,
which is at the intersection of a line---sometimes referred to as the crossbar---and 
a perpendicular ray, the  {\em stem}.
The group $\simGroup$ acts transitively on this \mabs, and as holotype we may take
the center $(0,0)$ and the lines $y=0$ and $\{x=0, y \le 0\}$.

Two such $\overline{\sf T}$'s can intersect in at most four points, 
 our standard analysis gives $R \le V_{\sC}+2n+1 \le 4\binom{n}{2} +2n+1$, and so
\beql{EqT1}
a_{\overline{\sf T}}(n) \le 2n^2 + 1~.
\eeq
Figure~\ref{FigT19} illustrates $a_{\overline{\sf T}}(3) = 19$.
However,  the configuration of three $\overline{\sf T}$'s shown in that figure is essentially unique,
and one can check (by hand, or by using Lemma~\ref{lemma1} below)  that there is no way to add a fourth $\overline{\sf T}$   that meets each of the first three
in four points. So equality does not hold in \eqref{EqT1} for $n \ge 4$. For $n=4$ the best we can 
achieve is 32 regions (Fig.\ \ref{FigT32}).

To treat the general case, we need to consider the directions of the stems. For two $\overline{\sf T}$'s to intersect in four points, their stems must not be too close to parallel. To see this, suppose the first $\overline{\sf T}$ is placed in some standard position---we take it to be $y=0$ and $\{x=0, y \le 0\}$ so that the stem is the negative $y$-axis. We claim that if another $\overline{\sf T}$ achieves the maximum number of intersections with this reference $\overline{\sf T}$, then its center node must lie in the lower half-plane $\{y < 0\}$. Indeed, if the center node of the second $\overline{\sf T}$ is $(x, y)$ with $y > 0$ and if the direction of its stem is $(\cos \theta, \sin \theta)$, then we must have $\sin \theta < 0$, or else the stem could not enter the lower half-plane. 
\begin{figure}[!htb]
\centering
\begin{tikzpicture}
	\rotatedcrosslines{0}{0}{30}
	\draw[dashed, thick] (0, 0) -- (3, 0);
	\node[scale=1.4] at (-1, 0.2) {$(x, y)$};
	\draw[fill=black] (0, 0) circle[radius=2pt];
	\draw (1, 0) arc[start angle = 0, end angle = 30, radius = 1];
	\node[scale=1.4] at (1.3, 0.4) {$\theta$};
\end{tikzpicture}
\caption{Parameters for a constrained long-legged $\overline{\sf T}$: coordinates $(x, y)$ of the center, and the angle $\theta$ between the stem and the $x$-axis.}
	\label{FigTpars}
\end{figure}

More precisely, if the center node lies in the first quadrant, both $\cos \theta$ and $\sin \theta$ must be negative. Assuming this and computing the point of intersection between the crossbar of the second $\overline{\sf T}$ and the $y$-axis, we find that this intersection point has $y>0$, so that the crossbar does not intersect the stem of the reference $\overline{\sf T}$. Similarly, assuming the second center node to lie in the second quadrant again forces the second crossbar to miss the first stem (if the center node lies on the positive $y$-axis and the stems intersect we would have a degenerate configuration, which we rule out). It follows that the second center node must lie in the lower half-plane, and then for the second stem to intersect the first crossbar (which is the $x$-axis), its direction $(\cos \theta, \sin \theta)$ must satisfy $\sin \theta > 0$. Remembering that our reference $\overline{\sf T}$ has $\theta = -\pi/2$, we have shown: 
\begin{lemma}
	If two $\overline{\sf T}$'s intersect in four points, and if their stems have the directions $(\cos \theta_1, \sin \theta_1)$ and $(\cos \theta_2, \sin \theta_2)$ respectively, then $\cos(\theta_1 - \theta_2) < 0 $. In other words, the angle between their bearings must be strictly greater  than $\pi/2$. \label{lemma1}
\end{lemma}
\begin{figure}[!htb]
	\centering
	\begin{tikzpicture}
		\rotatedcrosslines{-1.5}{0}{-45}
		\rotatedcrosslines{0}{0}{200}
	\end{tikzpicture}
	\caption{A good pair: two $\overline{\sf T}$'s that intersect in four points. The angle between the stems must be strictly greater that $\pi/2$.}
	\label{FigGP}
\end{figure}

This fact will determine how many pairs of $\overline{\sf T}$'s can achieve four intersections. Let us say that two $\overline{\sf T}$'s form a \emph{good pair} if they intersect in four points (Fig.\ \ref{FigGP}). We will denote the maximum number of good pairs among $n$ $\overline{\sf T}$'s by ${\sf gp}(n)$. We already know that with $n = 3$ we can have three good pairs, but that with $n = 4$, only five out of the six pairs can be good. Thus ${\sf gp}(3) = 3$ and ${\sf gp}(4) = 5$. 

To find a formula for ${\sf gp}(n)$, let us imagine that we choose $n$ points on the unit circle and form the graph
$\Omega$  whose vertices are those points and whose edges are the good pairs (points strictly more than $\pi/2$ apart). Then ${\sf gp}(4) = 5$ implies that $\Omega$ contains no $K_4$, that is, it contains no complete subgraph on four vertices. 

A classic theorem of Turán \cite[p. 108]{Bol98} then implies that $\Omega$ can have at most as many edges as the Turán graph  $T(3, n)$, the graph formed by letting the $n$ vertices be the set of integers $\{0, 1 \dots n-1\}$ 
and joining two vertices with an edge if and only if they are not congruent modulo $3$.  We can in fact achieve this bound simply by taking the $k$th angle $\theta_k$ to be $2\pi k/3$. Thus ${\sf gp}(n)$ equals the number of edges in the Turán graph, that is,
\begin{align*}
	{\sf gp}(3m) &= 3m^2,\\
	{\sf gp}(3m + 1) &= 3m^2 + 2m,\\
	{\sf gp}(3m + 2) &= 3m^2 + 4m + 1.
\end{align*}
It remains to check that we can turn this into a non-degenerate configuration, but this is easy: 
place the $k$th $\overline{\sf T}$  with its center vertex at $(\cos 2\pi k/3, \sin 2\pi k/3)$ and with its stem through the origin. Displace each group of $\overline{\sf T}$'s slightly so as to break up the triple intersection at the origin. Then spread the centers in each group so that each pair has one stem-crossbar intersection, as indicated in Figure~\ref{FigDisplace}. Finally, rotate the $\overline{\sf T}$'s by small generic angles to create the missing stem-stem and crossbar-crossbar intersections. 

\begin{figure}[!htb]
	\centering
	\begin{tikzpicture}
		\crosslines{0}{0}
		\crosslines{0.2}{0.2}
		\crosslines{0.4}{0.4}
	\end{tikzpicture}
	\caption{A step in breaking up the groups of overlaid $\overline{\sf T}$'s.}
	\label{FigDisplace}
\end{figure}

In the resulting configuration, two $\overline{\sf T}$'s in different groups have $4$ intersection points, while two $\overline{\sf T}$'s in the same group have $3$ intersection points, for a total of
$$
V_{\sC} = 4{\sf gp}(n) + 3\left(\binom{n}{2} - {\sf gp}(n)\right)
$$
intersections. Optimality is clear, since we have the maximum number of good pairs (contributing $4$ each), and all remaining pairs contribute as much as they can, namely $3$. Re-expressing the function ${\sf gp}$ using the floor function and simplifying the formula, we have shown that
\beql{EqT}
a_{\overline{\sf T}}(n) = \left\lfloor \frac{n^2}{3} \right \rfloor + \frac{3}{2}n^2 + \frac{1}{2}n + 1 = \left\lfloor \frac{11n^2 + 3n + 6}{6}\right\rfloor \quad (\seqnum{A389614})~.
\eeq

\begin{figure}[!htb]
 \centerline{\includegraphics[angle=0, width=3.00in]{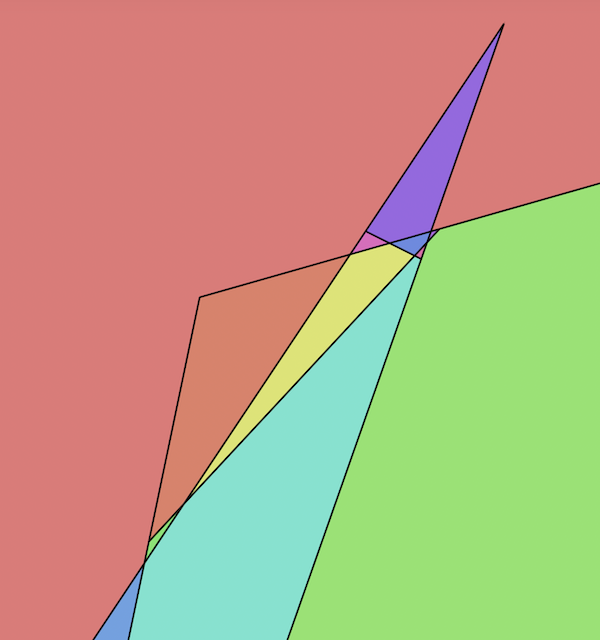}}
 \caption{Two constrained $\overline{\sf A}$'s  can  intersect in 8 points, and divide the plane into $a_{\overline{\sf A}}(2) = 13$ regions.
 (In this picture one $\overline{\sf A}$ has a purple tip, the other a brown tip.
 One of the intersection points is not shown: it occurs a long way below the bottom
 of the picture.  The arrowheads are also not shown.)}
 \label{FigsA2.13}
 \end{figure}

 \begin{figure}[!htb]
\centering
\begin{tikzpicture}[scale=.7, transform shape]   
\coordinate (a1) at (15.516576, 4.381652);
\coordinate (a2) at (11.58115260, 3.298787708);
\coordinate (a3) at (11.97140029, 2.358812006);
\coordinate (a4) at (0., 0.112137801);
\coordinate (a5) at (0., -4.471943158);
\coordinate (b1) at (19.452315, 3.885769);
\coordinate (b2) at (11.75354641, 3.429908786);
\coordinate (b3) at (11.92744506, 2.196041577);
\coordinate (b4) at (0., 2.733956684);
\coordinate (b5) at (0., -0.482294549);
\coordinate (c1) at (10.175527, -3.701396);
\coordinate (c2) at (12.15349946, 3.507635981);
\coordinate (c3) at (2.725865347, -3.080865250);
\coordinate (c4) at (12.83733987, 6.000000000);
\coordinate (c5) at (0., -2.853810204);

\draw[ -Stealth, thick] (a1) -- (a4);
\draw[ -Stealth, thick] (a1) -- (a5);
\draw[ thick] (a2) -- (a3);
\draw[ -Stealth, thick, red] (b1) -- (b4);
\draw[ -Stealth, thick, red] (b1) -- (b5);
\draw[ thick, red] (b2) -- (b3);
\draw[ -Stealth, thick, blue] (c1) -- (c4);
\draw[ -Stealth, thick, blue] (c1) -- (c5);
\draw[ thick, blue] (c2) -- (c3);
                \end{tikzpicture}
\caption{Three constrained $\overline{\sf A}$'s   
can divide the plane into 30 regions.}
 \label{FigcA3.30}
\end{figure}
 \begin{figure}[!htb]
\centering
\begin{tikzpicture}[scale=2.4, transform shape]   

\begin{scope}
\clip (10.5,2) rectangle (16,4.5);
\coordinate (a1) at (15.516576, 4.381652);
\coordinate (a2) at (11.58115260, 3.298787708);
\coordinate (a3) at (11.97140029, 2.358812006);
\coordinate (a4) at (0., 0.112137801);
\coordinate (a5) at (0., -4.471943158);
\coordinate (b1) at (19.452315, 3.885769);
\coordinate (b2) at (11.75354641, 3.429908786);
\coordinate (b3) at (11.92744506, 2.196041577);
\coordinate (b4) at (0., 2.733956684);
\coordinate (b5) at (0., -0.482294549);
\coordinate (c1) at (10.175527, -3.701396);
\coordinate (c2) at (12.15349946, 3.507635981);
\coordinate (c3) at (2.725865347, -3.080865250);
\coordinate (c4) at (12.83733987, 6.000000000);
\coordinate (c5) at (0., -2.853810204);

\draw[ -Stealth, thick] (a1) -- (a4);
\draw[ -Stealth, thick] (a1) -- (a5);
\draw[ thick] (a2) -- (a3);
\draw[ -Stealth, thick, red] (b1) -- (b4);
\draw[ -Stealth, thick, red] (b1) -- (b5);
\draw[ thick, red] (b2) -- (b3);
\draw[ -Stealth, thick, blue] (c1) -- (c4);
\draw[ -Stealth, thick, blue] (c1) -- (c5);
\draw[ thick, blue] (c2) -- (c3);
\end{scope}
                \end{tikzpicture}
\caption{Enlargement of the portion of Fig.\ \ref{FigcA3.30}
where the three crossbars come together.}
 \label{FigcA3.30crop}
\end{figure}

\subsection{The constrained long-legged \texorpdfstring{$\overline{{\sf A}}$}{{\sf A}}}\label{SeccA}
The {\em unconstrained}  long-legged ${\sf A}$ was discussed in \S\ref{SecAG}, where we showed that it is
equivalent to a three-armed {\sf V}.
A  {\em constrained long-legged} $\overline{\sf A}$  is defined
to be  a long-legged ${\sf A}$ in which the endpoints of the crossbar are equidistant from the tip.
This family is not holotypic: the group $\simGroup$ does not act transitively, and $\affineGroup$ does not preserve the constraint.

For both the unconstrained ${\sf A}$ and the  constrained $\overline{\sf A}$, the degrees of the base nodes 
are not all equal.
For both shapes, there are three base nodes: the tip, which has degree 2,
and the two crossbar nodes, which have degree 3. Our basic equation \eqref{EqEV2}
then gives
\beql{EqcA0}
R = V_{\sC} +2n + 1\,. 
\eeq

Two constrained $\overline{\sf A}$'s  can intersect in a maximum of $\kappa(\overline{\sf A}) = 8$ points,
and divide the plane into $R = a_{\overline{\sf A}}(2) = 13$ regions. To see this, observe that a line can cut an
$\overline{\sf A}$ in at most three points, so $\kappa(\overline{\sf A})  \le 9$. But if the two legs of a second  $\overline{\sf A}$  cut the first  $\overline{\sf A}$  in three points, it is impossible for the crossbar of the second $\overline{\sf A}$  to cut the first  $\overline{\sf A}$ in three points. So $\kappa(\overline{\sf A}) = 8$.  

Even the  construction that achieves  $\kappa(\overline{\sf A}) = 8$  was not easy to find (it was discovered by our program), so we provide a colored illustration
in Fig.\ \ref{FigsA2.13}.

From \eqref{EqcA0} we obtain
\beql{EqcA1}
a_{\overline{\sf A}}(n) \le 4n^2-2n+1 \quad (\seqnum{A054554}) ~.
\eeq
However, it appears that equality does not hold in \eqref{EqcA1} for $n > 2$.
We modified the program used in the previous section to apply to this
problem. For $0 \le n \le 5$ the maximum number of regions we found were
$1,3,13,30,53,83$.

In particular, for three $\overline{\sf A}$'s, only 23 crossings are possible, not 24, giving 30 regions rather than 31. We do not give full details of the proof of this claim, as they are elementary and not very enlightening, but in outline the argument goes as follows: suppose two $\overline{\sf A}$'s intersect in eight points. By enumerating all possible geometries of such a configuration, we can reduce them to two essentially different arrangements. 
An analysis of the geometry of these two cases gives bounds on the tip angles of the  $\overline{\sf A}$'s, and also on the angle between their symmetry axes. One then checks that these constraints are incompatible if we seek a configuration of three $\overline{\sf A}$'s, every pair of which intersect in eight points. Thus $24$ intersection points are impossible, and from examples we know that $23$ are possible.

This is as far as we have been able to determine the exact solutions to this problem. 
For $n = 0, 1, 2, 3$ the values of $a_{\overline{\sf A}}(n)$ are
\beql{EqcA2} 
1, 3, 13, 30 \quad (\seqnum{A397182}) ~.
\eeq
The OEIS entry \seqnum{A397182} has numerical coordinates for these arrangements, and will be
updated if and when further values are determined.
The \mabs\ $\overline{\sf A}$ was by far the most difficult that we studied.

We can, however, give reasonably close upper and lower bounds on $a_{\overline{\sf A}}(n)$.
We describe them for a general shape ${\sf S}$, since they may have more general applicability.

The lower bounds are obtained as follows. Suppose we have a configuration of $n$ copies of a shape ${\sf S}$, labeled $1$ through $n$, and that we assign integer multiplicities $m_1, \dots m_n$ to those figures. We think of such a figure-with-multiplicity as a cluster of $m_i$ superimposed copies of figure number $i$, and the whole collection of clusters as a configuration of $N = m_1 + \dots m_n$ figures. This is not yet a valid configuration, because of the superpositions,  but then we perturb the members of each cluster so as to remove the superpositions. 

In all the cases we are interested in, it is possible to spread the figures in each cluster so that each pair has $\xi({\sf S})$ intersection points, where $\xi({\sf S})$ is a  {\em local} intersection number defined as follows.
Recall that in \S\ref{SecPG} we defined the intersection number  $\kappa({\sf S})$ 
to be the maximum number of intersections between two copies of a figure ${\sf S}$. 
In those cases where the set of figures is generated from a holotype  by the action of a group, 
$\kappa({\sf S})$ is the maximum number of intersections between ${\sf S}$ and $g({\sf S})$ as $g$ ranges over the whole group; 
and we now define $\xi({\sf S})$ to be the corresponding 
maximum intersection number when $g$ is restricted to an arbitrarily small neighborhood of the identity element. 

Clearly $\xi({\sf S}) \leq \kappa({\sf S})$, and in general the inequality is strict. For example, the operation illustrated in Figure~\ref{FigDisplace} shows that while $\kappa(\overline{\sf T}) = 4$, we have $\xi(\overline{\sf T}) = 3$. Similarly $\kappa(\overline{\sf A}) = 8$ but $\xi(\overline{\sf A}) = 6$. 

To compute the total number of intersection points thus achieved, we first write down a table of the actual intersection counts between the figures ${\sf S}_1, \dots {\sf S}_n$,  before we assign any multiplicities. So as not to double-count, we only consider pairs $({\sf S}_i, {\sf S}_j)$ with $i<j$.   For example, the configuration with three $\overline{\sf A}$'s in Figure~\ref{FigcA3.30} produces Table~\ref{TableA9}.
\begin{table}[!htb]
	$$
	\begin{array}{l|llll}
		 & 1 & 2 & 3 \\
		\hline
		1 &  \ast & 8  & 7  \\
		2 &    &  \ast & 8\\
		3 &    &    & \ast \\
	\end{array}
	$$
	\caption{Intersections counts  between the 3 $\overline{\sf A}$'s in Figure~\ref{FigcA3.30}}\label{TableA9}
\end{table}
Assigning the figures the multiplicities $(m_1, m_2, m_3)$ and perturbing as described above, the total number of intersection points becomes
$$
I(m_1, m_2, m_3) = 8m_1m_2 + 8m_2m_3 + 7m_1m_3 + 6\left(\binom{m_1}{2} + \binom{m_2}{2} + \binom{m_3}{2}\right).
$$
By fixing $N = m_1 + m_2 + m_3$ and maximizing the quantity $I(m_1, m_2, m_3)$, we obtain lower bounds on the optimal intersection numbers, as follows:
$$
I(1, 2, 1) = 45, \text{ so } a_{\overline{\sf A}}(4) \geq 54,
$$
and 
$$
I(1, 2, 2) = 74, \text{ so } a_{\overline{\sf A}}(5) \geq 85,
$$
and so on. 

In general, then, we start with a graph formed from the intersection of ${\sf S}_1, \dots {\sf S}_n$, 
having multiplicities $m_i$, whose pairwise intersection counts are the coefficients of the terms $m_i m_j$, 
and the quantity $\xi$ determines how many additional intersections we get from each cluster. 
By fixing $N = m_1 + m_2 + \dots +  m_n$ and maximizing the corresponding quantity $I(m_1, m_2, \dots, m_n)$, we obtain lower bounds on the optimal intersection numbers.

Turning now to upper bounds, we make the following (trivial) observation: if we have a configuration and form its intersection table; and if we then choose a subset of the indices and delete the corresponding rows and columns, what remains is a possible intersection table. In other words, a principal submatrix of a possible intersection table is again a possible intersection table. This is obvious since the operation of dropping rows and corresponding columns is just the operation of omitting some of the figures in a configuration.  

As soon as we have an upper bound on the number of intersections for some $n$, this can be leveraged to obtain bounds for all larger values of $n$. We illustrate this idea by an example. Suppose we find a configuration of four $\overline{\sf A}$'s and write down its intersection table, as in Table~\ref{hypothetical}.
\begin{table}[!htb]
	$$
	\begin{array}{l|llll}
		& 1 & 2 & 3 & 4\\
		\hline
		1 &  \ast & a_{1, 2}  & a_{1, 3}  & a_{1, 4}\\
		2 &   &  \ast & a_{2, 3} & a_{2, 4}\\
		3 &   &   & \ast &  a_{3, 4}\\
		4 & & & & \ast
	\end{array}
	$$
	\caption{Hypothetical intersection table for four $\overline{\sf A}$'s.}\label{hypothetical}
\end{table}
Then, because dropping the first row and column  gives a valid intersection table for three  $\overline{\sf A}$'s, we have
$$
a_{2, 3} + a_{2, 4} +  a_{3, 4} \leq 23, 
$$
and for the same reason we have 
\begin{align*}
	a_{1, 3} + a_{1, 4} +  a_{3, 4} \leq 23, \\
	a_{1, 2} + a_{1, 4} +  a_{2, 4} \leq 23,\\
	a_{1, 2} + a_{1, 3} +  a_{2, 3} \leq 23.
\end{align*}
Adding up  these inequalities, we have
\beql{Eq4A}
2\sum_{i < j} a_{i, j} \leq 4\cdot 23.
\eeq

To formalize this, fix a shape ${\sf S}$, and let $M(n)$ denote the largest possible number of intersection points between $n$ copies of ${\sf S}$, that is, the largest possible sum of all entries in all possible intersection tables of size $n$ (we may assume that the  below-diagonal elements are zero). Thus, since Table~\ref{TableA9} is optimal, we have $M(3) = 23$ for the figure $\overline{\sf A}$, and our argument shows that this implies $M(4)\leq 46$. 

In general, we can imagine starting with a table of size $n$ and deleting one row and column pair in all possible ways,
obtaining  $n$ possible tables of size $n-1$. Next, we add the entries in each  table, and add up all these  sums. The result will be a multiple of the sum of the original table of size $n$, namely $n - 2$. This is so because any entry in the table, say $(1, 2)$, will be present precisely in those tables where the $(n-1)$-subset of remaining indices contains $1$ and $2$ (the diagonal elements are zero so the indices may be assumed to be distinct). On the other hand, each of the $n$ tables of size $n-1$ sums to at most $M(n-1)$. Since the original table was arbitrary (apart from being possible), we have
$$
(n - 2)M(n) \leq n M(n-1)
$$
as in \eqref{Eq4A}, and we conclude that
\beql{EqnA}
M(n) \leq \left \lfloor \frac{n}{n-2}M(n-1) \right \rfloor.
\eeq
(It would be possible to delete any number of columns,  but deleting a 
single column seems to give the best bound.)

For example, knowing that $M(3) = 23$  for $\overline{\sf A}$,  we get $M(4) \leq 46$, whence $M(5)\leq 76$ after taking the integer part, and so on. The upper and lower bounds are summarized in Table~\ref{TableAbounds}.

\begin{table}[!htb]
	$$
	\begin{array}{l|rr}
		n & \text{lower} & \text{upper}\\
		\hline
		1 & 3 & 3\\
		2 &  13 & 13\\
		3 &  30 & 30\\
		4 &  54 &  55\\
		5 & 85 & 87\\
		6 & 123 & 127\\
		7 & 169 & 174
	\end{array}
	$$
	\caption{Bounds on $a_{\overline{\sf A}}(n)$.}\label{TableAbounds}
\end{table}

It is also easy to extract bounds on the leading term of $M(n)$ as a function of $n$. For a lower bound, we maximize the function $I(m_1, m_2, m_3)$ over the real numbers, getting a leading term of $25n^2/7$. Rounding 
to an integer only incurs a linear correction.
For an upper bound, we may simply leave out the floor function from some point on --- $k$, say. 
Then the bound on $M(n)$ is
$$
M(n) \leq \frac{n(n-1)}{k(k-1)} M(k)
$$
for any fixed $k$. For example, using $M(7)\leq 159$ shows that the leading coefficient is at most $53/14 = 3.785\dots$. It is easy to improve this slightly by going past points where the floor would have made a difference. Moreover, the difference between the number of crossings and the number of regions is linear in $n$. Altogether, then, we find that $a_{\overline{\sf A}}(n)$ asymptotically satisfies
$$
\frac{25}{7}n^2 \leq a_{\overline{\sf A}}(n) < \frac{53}{14}n^2.
$$
It is conceivable that the lower bound is tight, but to prove this seems to require a closer analysis of the possible intersection tables.

\section{Finite shapes 1: Polygons}\label{SecFS}

\subsection{Convex polygons}\label{SecPolygon}
The \mabs\ of $k$-sided convex polygonal knives, for $k \ge 3$, abbreviated $k{\sf P}$, is the last of our \mabp\ in the affine genus.
Since $V_{\sB} = k$ and  $v_{\sB} = 2$, \eqref{EqEVf} tells us that $R =  V_{\sC} + 2$. Observe that any line segment can intersect a convex polygon in at most two points. In particular, each side in one of the polygons can intersect at most two sides of the other polygon. Therefore, two such polygons can intersect in at most $2k$ points, so $V_{\sC} \le 2k\binom{n}{2}$. 
This is easily achieved, so 
\beql{EqkP1}
a_{k{\sf P}}(n) = k n^2 - kn + 2 ~(n \ge 1), \mbox{~with~} a_{k{\sf P}}(0)=1~.
\eeq
The solutions for $k=3$, $4$, and $5$ (for the triangle, square or convex quadrilateral, and pentagon)
are given in sequences  \seqnum{A077588}, \seqnum{A069894}, and \seqnum{A386485}.
It was a slight surprise to discover that taking $k=1$ and $k=2$ 
in \eqref{EqkP1} matched \seqnum{A386480}
and \seqnum{A051890}, which are the numbers of regions attainable with $n$ circles
(see \S\ref{SecOO}) and $n$ ellipses, respectively.
More surprising, however, was the case $k=8$, which we discuss in the next
subsection.

\begin{figure}[!htb]
\centering
\begin{tikzpicture}[scale=.85, transform shape]

\node[draw, thick, regular polygon, regular polygon sides=8, minimum size=5cm] {}; 
\node[draw, thick, regular polygon, regular polygon sides=8, rotate=22.5, minimum size=5cm] {}; 
                \end{tikzpicture}
\caption{Two regular octagons divide the plane into 18 regions.}
 \label{FigOctagon2}
\end{figure}

\subsection{The octagon and concave quadrilateral}\label{SecOctagon}
For a convex octagon, the values of $a_{8{\sf P}}(n)$ 
are $1$, $2$, $18$, $50$, $98$, $162$, $242, \ldots$, that is,  $2(2n-1)^2$ for $n>0$.
Figure \ref{FigOctagon2}  illustrates $a_{8{\sf P}}(2) = 18$.
Remarkably, this sequence also matches   \seqnum{A077591},
whose definition is the maximum number of regions obtainable in the plane
by drawing $n$ {\em concave} quadrilaterals. 
The latter is illustrated in Figure~\ref{FigInky} for 
 $n=2$ (and which also has 18 regions). Since there is no obvious connection between the two dissections, 
 the result {\em may} be  just  a coincidence.

\begin{figure}[!htb]
\centering
\begin{tikzpicture}[scale=.85, transform shape]   

\draw[thick] (0,0) -- (4,1) -- (1,0) -- (4,-1) -- (0,0);
\draw[thick] (2.5,2.5) -- (3.3,-1.5) -- (2.5,1.5) -- (1.5,-1.5) -- (2.5,2.5);
                \end{tikzpicture}
\caption{Two concave quadrilaterals also divide the plane into 18 regions.
Is it just a concidence tht the maximum number of pieces that
can be obtained with $n$ (convex) octagonal cuts is the same as what can be
obtained with $n$ concave quadrilateral cuts? We do not know! }
 \label{FigInky}
\end{figure}

\begin{figure}[!htb]
	\centering
	\begin{tikzpicture}[scale=2]
		\coordinate (A) at (1, -0.1);
		\coordinate (B) at (72: 1);
		\coordinate (C) at (144: 1.4);
		\coordinate (D) at (195: 1.5);
		\coordinate (E) at (288: 1);
		\draw[thick, name path=AC] (A) -- (C);
		\draw[thick, name path=CE] (C) -- (E);
		\draw[thick, name path=EB] (E) -- (B);
		\draw[thick, name path=BD] (B) -- (D);
		\draw[thick, name path=DA] (D) -- (A);
		\coordinate (O) at (-1.2, -1);
		\coordinate (F) at (1.9, 1);
		\coordinate (G) at (2, 0);
		\coordinate (H) at (-0.5, 1);
		\draw[red, name path=OF] ($ (O)!-0.2cm!(F) $) -- (F);
		\draw[red, name path=OG] ($ (O)!-0.2cm!(G) $) -- (G);
		\draw[red, name path=OH] ($ (O)!-0.2cm!(H) $) -- (H);
		\path[name intersections={of=AC and OF, by=I1}];
		\fill[red] (I1) circle (2pt);
		\path[name intersections={of=AC and OH, by=I2}];
		\fill[red] (I2) circle (2pt);
		\path[name intersections={of=CE and OH, by=I3}];
		\fill[red] (I3) circle (2pt);
		\path[name intersections={of=CE and OG, by=I4}];
		\fill[red] (I4) circle (2pt);
		\path[name intersections={of=CE and OF, by=I5}];
		\fill[red] (I5) circle (2pt);
		\path[name intersections={of=EB and OG, by=I6}];
		\fill[red] (I6) circle (2pt);
		\path[name intersections={of=EB and OF, by=I7}];
		\fill[red] (I7) circle (2pt);
		\path[name intersections={of=DA and OF, by=I8}];
		\fill[red] (I8) circle (2pt);
		\path[name intersections={of=DA and OH, by=I9}];
		\fill[red] (I9) circle (2pt);
		\path[name intersections={of=BD and OH, by=I10}];
		\fill[red] (I10) circle (2pt);
	\end{tikzpicture}
	\caption{No line can intersect a pentagram in more than four points.}
	\label{FigurePentagramIntersections}
\end{figure}
\begin{figure}[!htb]
\centering
\begin{tikzpicture}[scale=1, transform shape]   

\def\radius{3}
    \def\sides{5}
    \def\skip{2} 

    \draw[line width=1pt] (90:\radius)
        \foreach \x in {1,...,\sides} {
            -- ({90 - 360 * \x * \skip / \sides}:\radius)
        } -- cycle;

\begin{scope}[rotate=24]
    \draw[line width=1pt] (90:\radius)
        \foreach \x in {1,...,\sides} {
            -- ({90 - 360 * \x * \skip / \sides}:\radius)
        } -- cycle;
\end{scope}

\begin{scope}[rotate=48]
    \draw[line width=1pt] (90:\radius)
        \foreach \x in {1,...,\sides} {
            -- ({90 - 360 * \x * \skip / \sides}:\radius)
        } -- cycle;
\end{scope}

\begin{scope}[rotate=48]
\draw[blue, ultra thick, fill] (0,3.0) circle [radius=.1];
\draw[blue, ultra thick, fill] (-.4,1.81) circle [radius=.1];
\end{scope}
\draw[blue, ultra thick, fill] (.15,.93) circle [radius=.1];

                \end{tikzpicture}
\caption{Three pentagrams can divide the plane into 77 regions.}
\label{FigPenta3}
\end{figure}

\subsection{The pentagram and  hexagram}\label{SecPentagram}
A {\em regular} pentagram or hexagram is based on 5 or 6 equally spaced
points around a circle. To attain the maximum number of regions
we will use regular figures, but our analysis applies equally to non-regular figures.
A single  pentagram (regular or non-regular)  has 5 base nodes, the exterior vertices, with degree $d_B=2$,  
and 5 self-intersections, in the interior. To compute~$\kappa$ for the pentagram, we can reason as follows: suppose we place a pentagram in the plane, and draw any line not incident with any of its vertices, as illustrated in Figure~\ref{FigurePentagramIntersections}. It is clear that the line will intersect a side of the pentagram if and only if the vertices joined by that side lie on opposite sides of the line. It follows from this that no line can intersect the pentagram in more than four points, and a fortiori that no line segment can intersect it in more than four points. Thus the five sides of one pentagram can intersect another pentagram in no more than $20$ points, so $\kappa \leq 20$. But this is also achievable by taking two pentagrams with their vertices on the same circle and placed so that their vertices interlace as we go around the circle. 

Since $\kappa = 20$, with $n$ pentagrams we have $V_{\sC} \le 5n + 20 \binom{n}{2}$, and
 our basic equation \eqref{EqEVf} then gives
 $R \le V_{\sC} + 2  \le 5n + 20 \binom{n}{2} + 2 = 10n^2-5n+2$ (\seqnum{A383466}),
 a bound which is easily achieved using regular pentagons based
 on $5n$ equally spaced points around a circle.
 The construction is illustrated for $n=3$  in Fig.\ \ref{FigPenta3}. 
 To verify the number of regions, note that each of 
 the  base nodes is a vertex of a triangular region (indicated by 
 the three blue dots  in Fig.\ \ref{FigPenta3}), containing $2n-1$ cells. We were led to investigate pentagrams because  an optimal 5-chain (\S\ref{SeckC})  necessarily contains a pentagram (see Figs.\ \ref{Fig345chains} and \ref{FigMy1}).

A hexagram initially contains two disjoint triangles and is not a connected graph, 
which would violate 
the assumption that our shapes are connected, 
but it {\em is} connected if we take the self-intersections into account. 
The analysis is similar to that of the pentagram, so we just state the result:
the maximum number of regions that can be attained in the plane by drawing
$n$ hexagrams is $2(6n^2-3n+1)$ (\seqnum{A386477}). 
 
\section{Finite shapes 2: Three curved shapes}\label{SecFS2}
Although our basic equation \eqref{EqEVf} would not apply to curved shapes with arbitrarily many twists,
it does apply to the three shapes considered in this section. These \mabp\  all belong
to the transitive similarity subgenus.

\begin{figure}[!htb]
\centering
\begin{tikzpicture}[scale=1.5, transform shape] 
\coordinate(A) at (-1,0);
\coordinate(C) at (-0.333,0);
\coordinate(D) at (0,0);
\coordinate(E) at (+0.333,0);
\coordinate(G) at (+1,0);
\draw[thick] (A) circle [radius=1.5];
\draw[thick] (C) circle [radius=1.5];
\draw[thick] (E) circle [radius=1.5];
\draw[thick] (G) circle [radius=1.5];
\coordinate(I1) at (-.667,1.462);
\coordinate(I6) at (+.667,1.462);
\coordinate(I2) at (-.333,1.344);
\coordinate(I5) at (+.333,1.344);
\coordinate(I3) at (0,1.462);
\coordinate(I4) at (0,1.118);
\coordinate(J1) at (-.667,-1.462);
\coordinate(J6) at (+.667,-1.462);
\coordinate(J2) at (-.333,-1.344);
\coordinate(J5) at (+.333,-1.344);
\coordinate(J3) at (0,-1.462);
\coordinate(J4) at (0,-1.118);
\draw[black, ultra thick, fill] (I1) circle [radius=.05];
\draw[black, ultra thick, fill] (I6) circle [radius=.05];
\draw[black, ultra thick, fill] (I2) circle [radius=.05];
\draw[black, ultra thick, fill] (I5) circle [radius=.05];
\draw[black, ultra thick, fill] (I3) circle [radius=.05];
\draw[black, ultra thick, fill] (I4) circle [radius=.05];
\draw[black, ultra thick, fill] (J1) circle [radius=.05];
\draw[black, ultra thick, fill] (J6) circle [radius=.05];
\draw[black, ultra thick, fill] (J2) circle [radius=.05];
\draw[black, ultra thick, fill] (J5) circle [radius=.05];
\draw[black, ultra thick, fill] (J3) circle [radius=.05];
\draw[black, ultra thick, fill] (J4) circle [radius=.05];
\coordinate(B1) at (-2.414,.5);
\coordinate(B4) at (-0.414,.5);
\coordinate(B2) at (-1.747,.5);
\coordinate(B3) at (-1.081,.5);
\draw[red, ultra thick, fill] (B1) circle [radius=.05];
\draw[red, ultra thick, fill] (B2) circle [radius=.05];
\draw[red, ultra thick, fill] (B3) circle [radius=.05];
\draw[red, ultra thick, fill] (B4) circle [radius=.05];
\draw[black,thin] (-3,0) -- (+3,0);
\draw[black, thick, fill] (A) circle [radius=.02];
\draw[black, thick, fill] (C) circle [radius=.02];
\draw[black, thick, fill] (D) circle [radius=.02];
\draw[black, thick, fill] (E) circle [radius=.02];
\draw[black, thick, fill] (G) circle [radius=.02];
\node[below, xshift=4pt] at (A) {$C_1$};
\node[below, xshift=1pt] at (C) {$P$};
\node[above] at (D) {$0$};
\node[below, xshift=-1pt] at (E) {$Q$};
\node[below, xshift=-4pt] at (G) {$C_2$};
                \end{tikzpicture}
\caption{Four circles can divide the plane into  14 regions.
There are four base nodes (red) and twelve crossing nodes (black).
The two outer circles have centers $C_1 = (-1,0)$ 
and $C_2 = (1,0)$,
and the two inner circles have centers $P = (-1/3,0)$ and $Q=(1/3,0)$.
All the circles have radius $3/2$.
The $x$-axis is for reference only, it is not part of the graph.}
 \label{Figcircle4d}
\end{figure}

\subsection{The circle}\label{SecOO}
It is a classical result that the maximum number of regions that can be obtained with $n$ circular cuts is
\beql{EqCirc1} 
a_{\sf O}(n) = n^2 - n + 2.
\eeq
This is already mentioned in Steiner's 1826 paper \cite{Ste26}, and other references can be found in \seqnum{A014206}.  We discuss the circle here because it is an interesting application of our method, and also because a similar argument will be used in the next section.

In order to apply our method, we formally define the circle shape to consist of a distinguished point, the base, and a circle of arbitrary positive radius through that point.  Note that without the base point, Euler's formula itself \eqref{Euler1} fails for a single circle, since there is a single edge, no vertices, and two regions, and $E-V+2 = 3$, not $2$. But, after adding a base node, we have $E-V+2 = 1-1+2 = 2$, as it should be.

It is clear that two circles can intersect in at most two points, so with $n$ circles, 
maintaining our standard convention that no other circle can pass through a base node,
we have $V_{\sB} = n$, $V_{\sC} \le 2 \binom{n}{2}$, and \eqref{EqEVf} implies 
that $R \le  2 \binom{n}{2} +2$.

If we can guarantee that every pair of circles meet in two points, then equality holds in the latter expression, and establishes \eqref{EqCirc1}. This is easily accomplished.  For $n \ge 2$ we draw two circles of radius $3/2$ 
centered respectively at $C_1 = (-1,0)$ and $C_2 = (1,0)$, and then draw $n-2$ further circles
of the same radius centered at points equally spaced along the line $C_1 - C_2$. 
Figure~\ref{Figcircle4d} illustrates the construction for $n=4$. The base nodes (shown in red), whose only function is to ensure that each circle contains at least one vertex, can be placed anywhere on the circles, except at the intersection points.  In Fig.~\ref{Figcircle4d} we have lined them up in a row.

\begin{figure}[!htb]
\centering
\begin{tikzpicture}[scale=.50, transform shape]   

\draw[thick] (0,0) circle [radius=3];
\draw[thick] (6,0) circle [radius=3];
\draw[thick] (3,-1.1875) circle [radius=3];
\draw[thick] (0,4.298) arc(-170:-10:3);
                \end{tikzpicture}
\caption{Two figure 8's can divide the plane into 12 regions}
 \label{Fig8.1}
\end{figure}

\begin{figure}[!htb]
	\centering
	\begin{tikzpicture}[scale=3]
		\pgfmathsetmacro{\theta}{35}
		\coordinate (A) at ({-1/2 - cos(\theta)/2}, {-sin(\theta)/2});
		\coordinate (B) at ({-1/2 + 3*cos(\theta)/2}, {3*sin(\theta)/2});
		\draw[thick] (-1, 0) circle[radius = 1];
		\draw[thick] (1, 0) circle[radius = 1];
		\draw[thick, dashed] (-1, 0) -- (1, 0);
		\draw[fill=black] (-1, 0) circle[radius=1pt];
		\draw[fill=black] (1, 0) circle[radius=1pt];
		\draw[thick] (A) circle[radius = 1];
		\draw[thick] (B) circle[radius = 1];
		\draw[thick, dashed] (A) -- (B);
		\draw[fill=black] (A) circle[radius=1pt];
		\draw[fill=black] (B) circle[radius=1pt];
		\draw[fill=red] (-1/2, 0) circle[radius=1pt];
	\end{tikzpicture}
	\caption{Eight intersections between two figure 8's can be achieved by an arbitrarily small rotation.}
	\label{Fig8.3}
\end{figure}

\subsection{The figure {\sf 8}}\label{Sec88}

A figure {\sf 8} shape consists of two tangential circles of equal but otherwise arbitrary radius.
A graph with $n$ {\sf 8}'s contains $2n$ circles, and so $a_{\sf O}(2n)$ from the previous section  is
an upper bound on $a_{\sf 8}(n)$. However, we cannot achieve 
all the intersections needed 
to achieve that bound, since the two constituent circles of each {\sf 8}
cannot properly intersect. The effect of this constraint is simply
to reduce the number of crossings by 2 for each figure {\sf 8}, 
and so the solution for the {\sf 8}
satisfies $a_{\sf 8}(n) \le 4n^2-3n+2$.

We can achieve equality in that expression using  figure {\sf 8}'s of equal size, by a construction similar to that of the previous section. 
We start from a figure {\sf 8} in a standard position, say the two circles $(x\pm 1)^2 + y^2 = 1$, with centers $C_1 = (-1,0)$, $C_2 = (1,0)$. Imagine adding another copy of this figure, slightly displaced. Then, for the right circle to intersect the original left circle in two points, its center must have $x < 1$, and be situated to the left of $C_2$.
By symmetry, the new left circle must have its center to the right of $C_1$. At the same time, the distance between the centers must be constant. This suggests rotating the whole figure about a point on the line $C_1 - C_2$. We cannot use the origin as the center of rotation, for
then the point of tangency between the circles would become a multiple point.
On the other hand, any other interior point on $C_1 - C_2$  will do, and any small rotation about that point will create eight intersection points. 
By placing $n-1$ equally-spaced figure {\sf 8}'s  along $C_1 - C_2$, we 
obtain a total of $n$ figure {\sf 8}'s, each pair of which intersect in eight points.
We conclude that\footnote{Eq.\ \eqref{Eq88} was independently discovered by the referee of the first version of this article.} 
\beql{Eq88}
a_{\sf 8}(n) = 4n^2-3n+2  \quad (\seqnum{A386486})~. 
\eeq

 \begin{figure}[!htb]
\centering
\begin{tikzpicture}[rotate=0] [scale=1]

\begin{scope} [scale=0.5]
\coordinate(A) at (3.2,3.35);
\coordinate(B) at (8.9,3.50);
\draw[thick] (A) circle [radius=3.02];
\draw[thick] (B) circle [radius=3.05];
\coordinate(C) at (5.88,2.03);
\coordinate(Cp2) at (10.45,-0.4);
\draw[-Stealth, thick, red] (C) -- (Cp2);
\coordinate(D) at (6.25,2.05);
\coordinate(Dp) at (1.35,-.55);
\draw[-Stealth, thick, red] (D) -- (Dp);
\end{scope}
                \end{tikzpicture}
\caption{Two lollipops can intersect in at most 7 points and can
 divide the plane into 10 regions.
The stems of the lollipops are colored red.}
 \label{FigLPOP2}
\end{figure}

\begin{figure}[!htb]
\centering
\begin{tikzpicture}[scale=1.20, transform shape]   

\begin{scope}
\clip (2,0.5) rectangle (10.5,8.75);
\coordinate(A) at (3.5,3.75);
\coordinate(B) at (9.25,3.75);
\coordinate(C) at (6.5,8.40);
\draw[thick] (A) circle [radius=3];
\draw[thick] (B) circle [radius=3];
\draw[thick] (C) circle [radius=3];
\coordinate(D) at (6.15,5.2);
\coordinate(Dp) at (12.00,8.20);
\coordinate(E) at (6.60,5.2);
\coordinate(Ep) at (0.75,8.25);
\coordinate(F) at (6.45,5.41);
\coordinate(Fp) at (6.43,0.75);
\draw[-Stealth, thick, red] (D) -- (10.45, 7.42);
\draw[-Stealth, thick, red] (E) -- (2.05,7.6);
\draw[-Stealth, thick, red] (F) -- (Fp);
\end{scope}
                \end{tikzpicture}
\caption{With each pair intersecting in 7 points, three lollipops can divide the plane into 25 regions.}
 \label{FigLPOP3}
\end{figure}

 \subsection{The lollipop or qoppa}\label{SecLP}
 Although it does not quite meet our conditions, since it is both curved and infinite,
 we cannot resist ending with  the  lollipop shape (\qoppaletter).
This is essentially  the same as the archaic Greek letter  qoppa (Unicode symbols U+03D8, U+03D9).
Formally, a lollipop or qoppa  is a circle together with a ray (the {\em stem}) emanating from its center
but with the portion of the ray inside the circle omitted.  Again $\simGroup$ acts transitively on this \mabs. 
Besides its intrinsic appeal, the lollipop is also a possible candidate for a long-legged version of 
both the letters {\sf P} and {\sf Q}. (The Latin letter {\sf Q} is in fact a descendant of the qoppa.)

Two lollipops can meet in a maximum of $\kappa = 7$ points (Fig.\ \ref{FigLPOP2}).
Equation \eqref{EqEV2} applies to this shape, and implies that the 
 maximum number of regions possible with $n$ lollipops is 
 \beql{EqLP1}
 R \le \frac{7n^2-5n+2}{2} \quad (\seqnum{A389608})~.
 \eeq 
This can be achieved  for $n \le 3$ (Fig.\ \ref{FigLPOP3}). For $n = 4$, however, 
we will show that the maximum number of intersection points is $40$, and the maximum number of regions is $45$. This can be achieved by taking an optimal configuration of three very large lollipops and inserting a fourth in the small central cavity.  
Figure~\ref{FigureQoppa40} shows an enlargement of the central portion of the configuration.
In this close-up, the three original lollipops appear as  $\overline{\sf T}$ shapes.

\begin{figure}[!htb]
	\centering
	\begin{tikzpicture}[scale = 4]
	\begin{scope}
		\useasboundingbox (-0.3, -0.3) rectangle (2, 1.3);
		\clip (-0.3, -0.3) rectangle (2, 1.3);
		\draw[thick] (0.45, 0.4) circle[radius = 0.55];
		\draw[thick] (1, 0.4) -- (2, 0.4);
		\rotatedcrosslines{1.15}{0.65}{-114};
		\rotatedcrosslines{1.04881168}{-0.05}{96}
		\rotatedcrosslines{0}{0}{0}
	\end{scope}
	\node[font=\large] at (-0.4, 1.35) {$1$};
	\node[font=\large] at (0.0, 1.4) {$2$};
	\node[font=\large] at (-0.3, -0.1) {$3$};
	\node[font=\large] at (2, 0.5) {$4$};
	\end{tikzpicture}
	\caption{Enlargement of the central portion of a configuration of four lollipops with 40 intersections. Three of the lollipops are so large they
		appear as $\overline{\sf T}$ shapes. It seems unlikely, but it is true that because  the lines numbered 1, 2, and 3 are actually circles, when the lines are extended, lines 1 and 2 will cross a second time, lines 1 and 3 will cross twice,
		and line 4 will cross line 1 a second time and will cross line 3 twice.}
	\label{FigureQoppa40}
\end{figure}

It turns out that the same kind of analysis that  we applied to the $\overline{\sf T}$ applies here, although we are not able to solve it completely. 
Recall that we define a \emph{good pair} to be a pair of $\overline{\sf T}$'s in which
the stem of each intersects the other in two points, and that this restricts the relative position of the two figures. For lollipops, we say that two of them form a good pair if they achieve \emph{three or four} stem-circle intersections, as in Figure~\ref{FigGPQ}. 
\begin{figure}[!htb]
	\centering
	\begin{tikzpicture}
		\qoppa{0}{0}{0}
		\qoppa{2.5}{0.5}{180}
		\qoppa{0}{-3}{0}
		\qoppa{1.7}{-2.1}{210}
	\end{tikzpicture}
	\caption{Good pairs of lollipops, with three or four stem-circle intersections. 
As in Fig.\ \ref{FigGP}, the angle between the stems must be strictly greater that $\pi/2$.} \label{FigGPQ}
\end{figure}

As in \S\ \ref{SeccT},  by choosing  coordinates for one of the lollipops and writing down equations for the stem of the other, we can establish an analogue of Lemma~\ref{lemma1}:
\begin{lemma}
	If two lollipops form a good pair, that is, if together they have three or four stem-circle intersections, then the angle between the directions of their stems must be strictly greater than $\pi/2$.\label{lemma2}
\end{lemma}
Moreover, one can show that in the case of three stem-circle intersections, it is not possible for the two stems to intersect, while  four stem-circle intersections are compatible with both one stem-stem and two circle-circle intersections. 
We have not included this fact as  part of the lemma as it will not be used below, but it could be relevant for future investigations.

Using this lemma, we can  show that our configuration of four lollipops with $40$ intersection points is optimal. Indeed, all pairs of circles contribute two intersections, and all pairs of stems contribute one intersection. The stem-circle intersections are given in Table~\ref{TableQoppa40}.

\begin{table}[!htb]
	$$
	\begin{array}{c|cccc}
		& S_1 & S_2 & S_3 & S_4 \\ \hline
		C_1    & \ast & 0  & 2  & 2 \\
		C_2    & 2  & \ast & 2 & 2  \\
		C_3    & 2  & 2 & \ast  & 2 \\
		C_4    & 2  & 2 & 2 & \ast  \\
	\end{array}
	$$
	\caption{Intersections between circles and stems in our $40$-intersection configuration of four lollipops.}\label{TableQoppa40}
\end{table}

The only number in the table that could be improved is the zero, which indicates that circle 1 and stem 2 do not intersect. But if that entry were equal to 1 or 2, we would have six good pairs, and we know that five is optimal for four directions. Thus this configuration is optimal, as claimed.

In summary, the values of  $a_{\qoppaletter}(n)$ for $n = 0, \ldots, 4$ are
\beql{EqQoppa2} 
1,  2, 10, 25, 45 \quad (\seqnum{A389624}) ~.
\eeq

Just as with the letters $\overline{\sf T}$ and $\overline{\sf A}$, we can leverage this result to obtain both lower and upper bounds for general $n$. For the lower bound, we assign multiplicities to the four lollipops in our best configuration and 
calculate  the number of intersection points after the configuration has been made generic by a small perturbation. This is where the lollipop proves more difficult than the $\overline{\sf T}$: while the maximum number of intersections between two lollipops is $\kappa = 7$, the local maximum (achievable by arbitrarily small motions) is only $\xi = 4$. Because of this large decrease, we are unable to show that the resulting configurations are optimal. At present they certainly are the best that we know---see Table~\ref{TableQbounds}. \seqnum{A389624} gives further details about our construction, including a drawing of and
coordinates for a  40-crossing, 45-region optimal arrangement of four lollipops.

\begin{table}[!htb]
	$$
	\begin{array}{l|rr}
		n & \text{lower} & \text{upper}\\
		\hline
		1 & 2 & 2\\
		2 &  10 & 10\\
		3 &  25 & 25\\
		4 &  45 &  45\\
		5 & 71 & 72\\
		6 & 104 & 106\\
		7 & 142 & 146\\
		8 & 186 & 193
	\end{array}
	$$
	\caption{Bounds on $a_{\qoppaletter} (n)$.}\label{TableQbounds}
\end{table}

\vspace*{+.1in}

The lollipop section ended at this point on May 3, 2026, when this article was resubmitted for review.
On May 29, 2026, Brady Haran and N.J.A.S. released a video \cite{HS26} presenting
the above results.  Several viewers of the video (Tavor Gochman, Siddhartha Mahajan \& Paras Chopra) obtained sharpened bounds for the lollipop, using assistance from LLMs.
Some more novel and geometrically original obstructions were thereafter developed by Matthias Paulsen, without use of AI assistance; these obstructions alone led to bounds that improved beyond the sharpenings discovered in the initial agentic searches. Moreover, this result of Paulsen catalyzed subsequent investigation: using AI assistance again, Siddhartha contributed a nontrivial combinatorial argument leveraging Paulsen's work. Finally, Siddhartha and Tavor later independently and concurrently arrived at the final solution, again having extensively utilized the help of LLMs, but with the work depending foundationally on human contribution, e.g. Paulsen's.

We mention two of these results. First, the lower bounds mentioned in the
previous paragraph have been shown to be optimal for all $n \le 17$ and $n = 19$.
The values of  $a_{\qoppaletter}(n)$ for $n = 0, \ldots, 17$ are
\beql{EqQoppa3}
1, 2, 10, 25, 45, 71, 104, 142, 186, 237, 294, 356, 425, 500, 580, 667, 761, 859,
\eeq
and the value for $n=19$ is 1076.  
Although some of the arguments  were initially found by AI assistants, 
these bounds now have a traditional (humanly-readable) proof  \cite{MaCh26, Pau26}. 

The second result, which at present is AI-derived and has not yet been confirmed, is an
optimal choice for the multiplicities needed to construct a solution for any $n \ge 4$ from 
the $n=4$ solution.  This was found independently by Gochman and by Maharajan \& Chopra
using AI assistants. If confirmed, this will completely solve the lollipop problem.
The result will be that $a_{\qoppaletter}(n)$ for $n \ge 0$ is equal to $F(n) + 2n^2 -n +1$, where 
\beql{EqQoppan1}
F(n) ~=~ \max_{ \substack{ m_1+m_2+m_3+m_4 = n\\
                             m_1, m_2, m_3, m_4 \ge 0}}   \left[ 3 \sum_{i < j} m_i m_j - 2 \min_{i<j} m_i m_j \right]\,.
\eeq
It is hoped that a full report of all the new results will be available soon.

\section{Open problems}\label{SecOP}
\begin{enumerate}
\item
Find optimal graphs for the two shapes that we
left unfinished: the constrained long-legged letter $\overline{\sf A}$ (\S\ref{SeccA}),
and the lollipop (\S\ref{SecLP}).

\item
{\em Enumerative Geometry.} We have already mentioned (in~\S\ref{SeckC}) the sequence 
\seqnum{A090338}, which gives the number of ways to draw 
$n$ lines in general position in the Euclidean  plane, and which is also
the number of ways to cut the plane into the maximum
number of pieces with $n$ straight cuts. The values are known for $n \le 9$.
If we omit ``in general position'', 
and just ask for the number of ways to draw $n$ lines in the plane,
there are more possibilities, and the values
(\seqnum{A241600}) are known only for $n \le 7$.
For $n=3$, for example,  there are four possibilities: three parallel lines, two parallel lines and a line cutting both, three lines through a point, and three lines in general position.
There are still more ways if we restrict the question to a disk (a pizza, say).
For $0 \le n \le 6$ it appears that the values
are $1, 1, 2, 5, 19, 130, 1814$ (\seqnum{A272906}), although the last term is unconfirmed.
There has been no progress on
either of the last two questions for ten years.

\item
More generally, for any shape $S$ for which we have determined the maximum number 
of regions $a_S(n)$ (which includes all the entries in Tables~\ref{Table1} and \ref{Table2}),
we may ask for the number of distinct graphs having that many regions. 
The case when $S$ is a circle has been studied by J. Wild in \seqnum{A250001} and \seqnum{A288554}.
Wild shows that the numbers of distinct 
planar graphs  composed of $n$ overlapping circles that  divide the plane into $n^2-n+2$ regions
are respectively $1, 1, 4, 45, 5102$ for $n=1, \ldots, 5$. The four possibilities for $n=3$ are shown in
Fig.\ \ref{3circles4ways}. One of the $45$ is shown in Fig.\ \ref{Figcircle4d}.  Arrangements of circles in the plane 
are classified from a different point of view in Mathar \cite{Mat16}.

 \begin{figure}[!htb]
\centering
\begin{tikzpicture}[scale=.8, transform shape]   

\def\s1{-6.5}
\draw[thick] (\s1+1.5,6) circle [radius=.9];
\draw[thick] (\s1+2.35,6) circle [radius=.9];
\draw[thick] (\s1+1.85,5.3) circle [radius=.9];
\def\c1{.25}
\draw[thick] (\s1+4.5,\c1+6) circle [radius=.8];
\draw[thick] (\s1+6,\c1+6) circle [radius=.8];
\draw[thick] (\s1+5.21,\c1+4.7) circle [radius=.8];

\def\s2{3.4}
\draw[thick] (1.85,\s2+3) circle [radius=.95];
\draw[thick] (1.85,\s2+1.3) circle [radius=.95];
\draw[thick] (1.85,\s2+2.15) circle [radius=.65];
\def\c2{.3}
\draw[thick] (4.25,\s2+\c2+2.4) circle [radius=.85];
\draw[thick] (4.25,\s2+\c2+1.15) circle [radius=.85];
\draw[thick] (4.25,\s2+\c2+1.8) circle [radius=.35];

                \end{tikzpicture}
\caption{Three circles can divide the plane into a maximum of 8 regions, and this can be done in 4 ways. (Only the first is a Venn diagram.)}
\label{3circles4ways}
\end{figure}
 
 \item
 Is there a geometric explanation 
 for the apparent  coincidence mentioned in~\S\ref{SecOctagon}?
 
\item
{\em Geometrical Configurations.} In the classical theory of
Geometrical Configurations \cite{Gru09},  \cite[Chap.\  III, {\em Projective Configurations}]{HCV52}, \cite{Lev29},
a configuration with symbol $(p_{\lambda},  \ell_{\pi})$
is an arrangement of $p$ points and $\ell$ lines in the plane such that each point
is incident with $\lambda$ lines, and each line  is incident with $\pi$ points.

 Our graphs $\Gamma_{\sf S}(n)$ satisfy this condition provided the number of 
 crossings per arm ($\cpa$) is constant, as it is in most cases  considered here, 
 and if so then $\Gamma_{\sf S}(n)$ 
 is a geometrical configuration with $p = V_{\sC}$, $\lambda = 2$,
 $\ell =$ number of arms $\alpha(S)$, and $\pi = \cpa$.
 
 For example, the optimal constrained long-legged
 $\overline{\sf T}$  graphs discussed in~\S\ref{SeccT},
 which we showed do not exist for $n>3$,  would have formed
  configurations with $p = 4\binom{n}{2}$, $\lambda=2$, $\ell = 2n$, and $\pi=2(n-1)$
  for integers $n \ge 1$. Could the classical theory have resolved that question?
  Does the theory 
  have any other application to our problems?  Gr\"{u}nbaum's book \cite{Gru09},
  although it is mostly concerned with configurations of type $(p_{\lambda}, p_{\lambda})$,
  contains a large number of striking illustrations  similar to those in the present article.

\item {\em The Inverse Problem.} Can one characterize the  integer sequences that arise 
as $\{a_{S}(n): n \ge 0\}$ for some shape $S$?  They are obviously monotonically increasing, and 
\eqref{EqVx} gives a bound on the growth rate.

\item {\em Computational Complexity.} All our shapes have been semi-algebraic curves. This means that for any shape $S$ and for any integers $n$ and $k$, the statement ``it is possible to place $n$ copies of shape $S$ so as to create $k$ intersections'' is decidable \cite{Bas18}. Of course, general-purpose algorithms for quantifier elimination have doubly exponential complexity, so a frontal attack is futile. However, it is conceivable that certain classes of shapes lead to more tractable problems. After all, the shapes we are interested in lead to polynomial growth in $n$, while more general classes can give rise to exponentially many regions.  
\end{enumerate}

\section{Acknowledgments}
N.J.A.S. mentioned the subject of this article in a  Research Experience for Undergraduates
seminar at Rutgers University on June 24 2025.  Following the talk, 
three participants in the seminar,
Edward Xiong, Jonathan  Pei, and David O.\ H.\ Cutler,
solved the 3-armed {\sf V} problem (see \S\ref{SeckV}).
Their solution led us to  the edge-vertex formula of~\eqref{EqEV},
which underlies all our calculations. We thank
Susanna S.\ Cuyler,  Sean A.\ Irvine, Gareth McCaughan,  Pete McPartian, 
Hilarie Orman \& Rich Schroeppel, Scott R.\ Shannon, and Jonathan Wild for comments on the manuscript.
We also thank
Scott R.\ Shannon for his  dramatic illustrations
of the pentagram  and hexagram  graphs (see \seqnum{A383466} and  \seqnum{A386477}).

We are very grateful to the referee of the first version of this paper for many helpful comments and suggestions.

Finally, we thank Brady Haran and his Numberphile colleagues for creating the video \cite{HS26} which it seems will soon lead to a complete solution of the lollipop problem (\S\ref{SecLP}).


\bigskip
\hrule
\bigskip

\noindent 2020 {\it Mathematics Subject Classification}:
Primary 05C10; Secondary 05C63, 52C30.

\noindent \emph{Keywords: } pancake cutting, plane dissection, bent line, zig-zag, configuration, cookie cutter

\bigskip
\hrule
\bigskip

\noindent (Concerned with sequences
\seqnum{A000124},
\seqnum{A000125},
\seqnum{A000326},
\seqnum{A008865},
\seqnum{A014206},
\seqnum{A034856},
\seqnum{A046127},
\seqnum{A051890},
\seqnum{A054554},
\seqnum{A058331},
\seqnum{A069894},
\seqnum{A077588},
\seqnum{A077591},
\seqnum{A080856},
\seqnum{A084849},
\seqnum{A090338},
\seqnum{A117625},
\seqnum{A125201},
\seqnum{A130883},
\seqnum{A140063},
\seqnum{A140064},
\seqnum{A143689},
\seqnum{A152948},
\seqnum{A241600},
\seqnum{A250001},
\seqnum{A272906},
\seqnum{A288554},
\seqnum{A383464},
\seqnum{A383465},
\seqnum{A383466},
\seqnum{A386477},
\seqnum{A386478},
\seqnum{A386479},
\seqnum{A386480},
\seqnum{A386481},
\seqnum{A386485},
\seqnum{A386486},
\seqnum{A387525},
\seqnum{A389608},
\seqnum{A389614},
\seqnum{A389624},
\seqnum{A393441},
\seqnum{A393442},
\seqnum{A393448},
\seqnum{A397182}.)

\bigskip
\hrule
\bigskip

\vspace*{+.1in}
\noindent
Received
revised version received
Published in {\it Journal of Integer Sequences},

\bigskip
\hrule
\bigskip

\noindent
Return to \href{https://cs.uwaterloo.ca/journals/JIS/}{Journal of Integer Sequences home page}.
\vskip .1in


\begin{thebibliography}{99}

\bibitem{Ale81}
G. L. Alexanderson and J. E. Wetzel, 
Arrangements of planes in space, 
{\em  Discrete Math.} \textbf{34} (1981), 219--240.

\bibitem{Bar19}
J.-L. Baril  and C. Moreira Dos Santos,
 The pizza-cutter's problem and hamiltonian paths,
 {\em Math. Mag.} \textbf{92} (2019), 359--367.

\bibitem{Bas18}
S. Basu and B. Mishra,
Computational and quantitative real algebraic geometry,
in {\em Handbook of Discrete and Computational Geometry},
eds.  J. E. Goodman, J. O'Rourke and C. D. T\'{o}th,
Chapman and Hall/CRC Press, pp.  969--1002, 2018.

\bibitem{Bol98}  
B. Bollob\'{a}s, 
{\em Modern Graph Theory},
 Springer, 1998.

\bibitem{Com74} 
L. Comtet, {\em Advanced Combinatorics}, Reidel, 1974.

\bibitem{CDFP38}
H.~S.~M.~Coxeter, P.~du~Val, H.~T.~Flather, and J.~F.~Petrie,
{\em The Fifty-Nine Icosahedra},
University of Toronto Studies (Math. Series, No. 6), Toronto, 1938.
Reprinted, Springer, 1982.

\bibitem{Cro97}
P. R. Cromwell,
{\em Polyhedra}, Cambridge, 1997.

\bibitem{Ede86}
H. Edelsbrunner, J. O'Rourke, and R. Seidel,
Constructing arrangements of lines and hyperplanes with applications,
{\em SIAM J. Comput.} \textbf{15} (1986), 341--363.


\bibitem{Gar66}  
M. Gardner, 
{\em New Mathematical Diversions from Scientific American},
 Simon and Schuster, 1966.

\bibitem{GKP94}  
R. L. Graham, D. E. Knuth, and O. Patashnik, 
{\em Concrete Mathematics},
 Addison-Wesley,  2nd ed., 1994. 
 
 \bibitem{Gru67}  
B. Gr\"{u}nbaum, 
{\em Convex Polytopes},
John Wiley \& Sons, 1967.

\bibitem{Gru72}  
B. Gr\"{u}nbaum, 
{\em Arrangements and Spreads}, 
CBMS Regional Conference Series in Mathematics, No. 10,
Amer. Math. Soc., 1972.

\bibitem{Gru75}  
B. Gr\"{u}nbaum, 
Venn diagrams and independent families of sets,
{\em Math. Mag.} \textbf{48} (1975), 12--23.

\bibitem{Gru09}  
B. Gr\"{u}nbaum, 
{\em Configurations of Points and Lines}, 
Graduate Studies in Math. \textbf{103}, Amer. Math. Soc., 2009.

\bibitem{HS26}
B. Haran and N. J. A. Sloane, 
{\em The Unsolved Lollipop Problem,}
Youtube Numberphile video, June 2026,
\url{https://www.youtube.com/watch?v=v8e-tYey7ts}.

\bibitem{HCV52} 
D. Hilbert and S. Cohn-Vossen, {\em Geometry and the Imagination},
Chelsea, NY, 1952; 2nd ed., Springer, 1996. 

\bibitem{Kin09}
D. Kinsella,
{\em Plane division by Lines AND Circles (Problem, Analysis and Solution)},
A Maths Puzzle Page, no date, but modified in 2009. This document is
now apparently lost,
although an incomplete copy was preserved by the Internet Archive's
{\em Wayback Machine}
and can be seen here:
\url{https://web.archive.org/web/20210414225350/http://www.90thkilmacudscouts.com:80/maths/circles_lines_soln.html
}.

\bibitem{Lak76} 
I. Lakatos, 
{\em Proofs and Refutations: The Logic of Mathematical Discovery},
Cambridge, 1976.

\bibitem{Lev29} 
F. Levi, 
{\em Geometrische Konfigurationen}, Hirzel, 1929.

\bibitem{MaCh26}
S. Mahajan and P. Chopra, 
A two-graph refinement of Paulsen's lollipop bounds, arxiv preprint arXiv:2606.06064 
[math.CO], 2026. Available at \url{https://arxiv.org/abs/2606.06064}.

\bibitem{Mat16}
R. J. Mathar, Topologically distinct sets of non-intersecting circles in the plane,
arxiv preprint arXiv:1603.00077 [math.CO],  2016. Available at
\url{https://arxiv.org/abs/1603.00077}.

\bibitem{Mat02}
J. Matou\v{s}ek,
{\em Lectures on Discrete Geometry},
Springer, 2002.

\bibitem{Mit09}
T. S. Michael,
{\em How to Guard an Art Gallery and Other Discrete Mathematical Adventures},
Johns Hopkins University Press, 2009.

\bibitem{Mil82}
R. E. Miles, A generalization of a formula of Steiner,
{\em Z. f\"{u}r Wahrscheinlichkeitstheorie und Verwandte Gebiete}
\textbf{61} (1982), 375--378.

\bibitem{OEIS}  
OEIS Foundation Inc.\ (2025), 
{\em The On-Line Encyclopedia of Integer Sequences}, 
\url{https://oeis.org}. 

\bibitem{Pac95}
J. Pach and P. K. Agarwal,
{Combinatorial Geometry},
John Wiley \& Sons, 1995.

\bibitem{Pau26}
M. Paulsen, The lollipop problem, preprint, June 4, 2026,
\url{https://math.mp42.de/lollipop/}. 

\bibitem{Pol54}
G. P\'{o}lya,
{\em Induction and Analogy in Mathematics},
Princeton University Press, 1954.

\bibitem{Rin56}
G. Ringel, Teilungen der Ebene durch Geraden oder topologische Geraden,
{\em Math. Z.} \textbf{64} (1956)., 79--102.

\bibitem{RuWe05} 
F. Ruskey and M. Weston,
A survey of Venn diagrams,
{\em Electron. J. Combin.}
Dynamic Surveys \textbf{5}, 2005.

\bibitem{Ste26}
J. Steiner, Einige Gesetze \"{u}ber die Theilung der Ebene und des Raumes, 
{\em J. Reine Angew. Math.} \textbf{1} (1826), 349--364;
{\em Gesammelte Werke}, ed.  K. Weierstrass, Vol. 1,  pp. 80--100,  AMS Chelsea Publishing, 1997.

\bibitem{Wen71}
M. J. Wenninger, 
{\em Polyhedron Models},
Cambridge, 1971.

\bibitem{Wet78}
J. E. Wetzel,
 On the division of the plane by lines,
 {\em Amer. Math. Monthly}, \textbf{85} (1978), 647--656.
 
\bibitem{Yag87} 
A. M. Yaglom and I. M. Yaglom, 
{\em Challenging Mathematical Problems with Elementary Solutions. Vol. I. Combinatorial Analysis and Probability Theory}, 
Dover Publications, 1987.

\bibitem{Yea56}  
W. B. Yeats, 
{\em The Collected Poems of W. B. Yeats}, Macmillan, 1956.

\end{thebibliography}
\end{document}